\documentclass[12pt]{amsart}
\usepackage{amsmath,amssymb,amscd,latexsym,eufrak}

\newcommand{\mvskip}{\vspace{5mm}}
\newcommand{\svskip}{\vspace{3mm}}

\newcommand{\C}{{\Bbb C}}

\newcommand{\Z}{{\Bbb Z}}
\newcommand{\BP}{{\Bbb P}}
\newcommand{\Q}{{\Bbb Q}}

\newcommand{\Supp}{{\rm Supp}\:}
\newcommand{\Sing}{{\rm Sing}}

\newcommand{\ML}{{\rm ML}}

\newcommand{\red}{{\rm red}\:}

\newcommand{\QED}{{\unskip\nobreak\hfil\penalty50\quad\null\nobreak\hfil
{$\Box$}\parfillskip0pt\finalhyphendemerits0\par\medskip}}
\newcommand{\Proof}{\noindent{\bf Proof.}\quad}
\newcommand{\A}{{\Bbb A}}
\newcommand{\Spec}{{\rm Spec}\:}
\newcommand{\CSpec}{{\mathcal Spec}\:}
\newcommand{\Proj}{{\rm Proj}\:}
\newcommand{\lto}{\longrightarrow}

\newcommand{\lkd}{\ol{\kappa}}

\newcommand{\codim}{{\rm codim}\;}

\newcommand{\Pic}{{\rm Pic}\:}

\newcommand{\GL}{{\rm GL}}

\newcommand{\Ker}{{\rm Ker}\:}

\newcommand{\im}{{\rm Im}\:}
\newcommand{\id}{{\rm id}\:}

\newcommand{\ST}{{\mathcal T}}
\newcommand{\SL}{{\mathcal L}}

\newcommand{\SC}{{\mathcal C}}

\newcommand{\SO}{{\mathcal O}}
\newcommand{\SU}{{\mathcal U}}
\newcommand{\SV}{{\mathcal V}}

\newcommand{\gm}{\EuFrak{m}}

\newcommand{\gp}{\EuFrak{p}}

\newcommand{\GM}{\EuFrak{M}}

\newcommand{\Hom}{{\rm Hom}}

\newcommand{\Ext}{{\rm Ext}}
\newcommand{\st}[1]{\stackrel{{#1}}{\longrightarrow}}

\newcommand{\wt}{\widetilde}
\newcommand{\ol}{\overline}
\newcommand{\wh}{\widehat}
\newcommand{\quot}{/\!/}

\newtheorem{thm}{Theorem}[section]
\newtheorem{lem}[thm]{Lemma}
\newtheorem{prop}[thm]{Proposition}
\newtheorem{cor}[thm]{Corollary}

\newtheorem{remark}[thm]{{\sc Remark}}

\newtheorem{example}[thm]{{\sc Example}}

\newtheorem{question}[thm]{Question}

\begin{document}
\title[$\A^1_*$-fibrations on affine threefolds]{$\A^1_*$-fibrations on affine threefolds}
\author{R.V. Gurjar, M. Koras, K, Masuda, \\  M. Miyanishi and P. Russell}
\thanks{The second author was supported by a Polish grant NCN N201 608640. The third and fourth authors were 
respectively supported by Grant-in-Aid for Scientific Research (C) 22540059 and 21540055, JSPS. The fifth 
author was supported by a grant from NSERC, Canada}
%\date{December 6, 2011}
\maketitle

\begin{abstract}
The main theme of the present article is $\A^1_*$-fibrations defined on affine threefolds. The difference between 
$\A^1_*$-fibration and the quotient morphism by a $G_m$-action is more essential than in the case of an 
$\A^1$-fibration and the quotient morphism by a $G_a$-action. We consider necessary (and partly sufficient) 
conditions under which a given $\A^1_*$-fibration becomes the quotient morphism by a $G_m$-action. 
Then we consider flat $\A^1_*$-fibrations which are expected to be surjective, but this turns out to be not 
the case by an example of Winkelmann \cite{W}. This example gives also a quasi-finite endomorphism of $\A^2$ 
which is not surjective \cite{Jelonek}. We consider then the structure of a smooth affine threefold which has 
a flat $\A^1$-fibration or a flat $\A^1_*$-fibration. More precisely, we consider affine threefolds 
with the additional condition that they are contractible. 
\end{abstract}

\section*{Introduction}

In affine algebraic geometry, our knowledge on affine algebraic surfaces is fairly rich with various methods 
of studying them. Meanwhile, knowledge of affine threefolds is very limited. It is 
partly because strong geometric approaches are not available or still under development. 

A possible geometric approach is to limit ourselves to the case where the affine threefolds in consideration 
have fibrations by surfaces or curves, say the affine plane $\A^2$ or the affine line $\A^1$ or have group actions 
which give the quotient morphisms. To be more concrete, $\A^2$-fibrations were observed to give characterizations of 
the affine $3$-space $\A^3$ (see \cite{Mi4, Mi5, KZ2}). Meanwhile, the quotient morphism $q : Y \to Y\quot G_a$ 
for an affine threefold has been considered by many people including P. Bonnet, Sh. Kaliman, D. Finston,  
J.K. Deveney et al. See \cite{GMM} for the references. 

In \cite{GMM}, some of the present authors considered $\A^1$-fibrations and quotient morphisms for an action of $G_a$. 
There, a main point is that an $\A^1$-fibration is factored by a quotient morphism by $G_a$ and the study is 
essentially reduced to the quotient morphism by $G_a$. In this article, we consider $\A^1_*$-fibrations and  
quotient morphisms by $G_m$. Contrary to the case of $\A^1$-fibrations, $\A^1_*$-fibrations are not necessarily 
factored by a quotient morphism by $G_m$. The possibility to factor depends on the nature of the singular fibers of the 
$\A^1_*$-fibration (see Theorem \ref{Theorem 2.20}). The quotient morphism by $G_m$ provides us a fertile ground 
for research, and the $\A^1_*$-fibrations seem to be even more fertile, but mysterious. For example, singular fibers 
of $\A^1_*$-fibration are fully described in the surface case (see Lemma \ref{Lemma 2.2}) but not completely in the 
case of threefolds, and the locus of singular fibers is shown to be a closed set with the unmixedness condition on 
singular fibers (cf. Lemma \ref{Lemma 2.11}).

Our primary purpose in this article is to develop a theory of dealing with $\A^1_*$-fibrations by combining 
algebraic geometry with algebraic topology and commutative algebra and to apply it to elucidate the structure of
such objects as homology (or contractible) threefolds. Thus the article contains known results as well as 
what we hope to be original ones. 

Here we comment on the notation $\A^1_*$ for which the notation $\C^*$ is in more common use. But $\C^*$ is 
used to mean, for example, that an affine scheme $\Spec A$ has only constant invertible functions, i.e., $A^*
=\C^*$. It is better to distinguish the curve $\C^*$ from the multiplicative group $\C^*$. Hence we use $\A^1_*$ 
for the curve and $G_m$ for the group.  

This article is a product of the discussions which we had at McGill University in August, 2011 at the occasion 
of the Workshop on Complex-Analytic and Algebraic Trends in the Geometry of Varieties sponsored by the Centre de 
Recherches Math\'ematiques, Montr\'eal and held at Universit\'e de Quebec \`a Montr\'eal (UQAM). We are very grateful to 
the Department of Mathematics and Statistics, McGill University, and Professor Steven Lu of UQAM for providing us 
this nice opportunity for discussions. 

Finally the authors appreciate various comments by the referees. As pointed out by one of them, $\Q$-homology threefold 
might as well be called $\Q$-homology $3$-space. 

\section{Preliminaries}
We summarize various results which we make use of in the subsequent arguments. The ground field is 
always the complex number field $\C$. A {\em homology $n$-fold} is a smooth affine algebraic variety $X$ of 
dimension $n$ such that $H_i(X;\Z)=0$ for every $i > 0$. If $H_i(X;\Q)=0$ for every $i > 0$ instead, we call 
$X$ a {\em $\Q$-homology $n$-fold}. 

\begin{lem}\label{Lemma 1.1}
Let $X=\Spec A$ be a homology $n$-fold. Let $(V,D)$ be a pair of a smooth projective variety $V$ and 
a reduced effective divisor $D$ on $V$ such that $X$ is a Zariski open set of $V$ with $D=V-X$ and 
$D$ is a divisor with simple normal crossings. Then the following assertions hold.
\begin{enumerate}
\item[(1)]
$X$ is factorial and $A^*=\C^*$. 
\item[(2)]
$H^1(V,\SO_V)=H^2(V,\SO_V)=0$.
\end{enumerate}
\end{lem}
\Proof
Similar considerations can be found in \cite{Fujita}. To simplify the arguments, we treat the case $n=3$. 
Consider the exact sequence of $\Z$-cohomology groups for the pair $(V,D)$,
\[
\begin{array}{lllllllll}
0 &\to& H^0(V,D) &\to& H^0(V) &\to& H^0(D) &\to& H^1(V,D) \\
  &\to& H^1(V) &\to& H^1(D) &\to& H^2(V,D) &\to& H^2(V) \\
  &\to& H^2(D) &\to& H^3(V,D) &\to& H^3(V) &\to& H^3(D) \\
  &\to& H^4(V,D) &\to& H^4(V) &\to& H^4(D) &\to& H^5(V,D) \\
  &\to& H^5(V) &\to& H^5(D) 
\end{array}
\]
By Lefschetz duality, we have $H^i(V,D) \cong H_{6-i}(X)$ for $0 \le i \le 6$. Since $X$ is a homology manifold, 
$H_i(X)=0$ for $1 \le i \le 3$ and since $X$ is affine, $H_i(X)=0$ for $4 \le i \le 6$ 
(see \cite[Theorem 7.1]{Milnor}). Hence we obtain the isomorphism $H^i(V) \cong H^i(D)$ for $0 \le i < 6$. 
Since $H^5(D)=0$ as $\dim D=2$, we have $H^5(V;\Z)=0$. By Poincar\'e duality, we have $H_1(V;\Z)=0$.  
Then the universal coefficient theorem \cite[Theorem 3, p.243]{Spanier} implies that $H^1(V;\Z)=0$, and hence 
$H^1(V;\C)=0$. The Hodge decomposition then implies that $H^1(V,\SO_V)=0$. Now consider an exact sequence
\[
0 \lto \Z \lto \SO_V \st{\exp(2\pi i\ )} \SO_V^* \lto 0
\]
and the associated exact sequence 
\begin{equation}
H^1(V,\SO_V) \to H^1(V,\SO_V^*) \to H^2(V;\Z) \to H^2(V,\SO_V)
\end{equation}

On the other hand, consider the isomorphism $H^4(V;\Z) \cong H^4(D;\Z)$. By the Mayer-Vietoris exact sequence, 
it follows that $H^4(D;\Z)$ is a free abelian group generated by the classes of the irreducible components of 
$D$. The Poincar\'e duality and the universal coefficients theorem implies
\[
H^4(V;\Z) \cong H_2(V;\Z) \cong H^2(V;\Z),
\]
where we note that $H_1(V;\Z)=0$. Hence we know that $H^2(V;\Z)$ is a free abelian group generated by the 
first Chern classes of the irreducible components of $D$. Since all these classes are algebraic, we conclude 
that $H^2(V;\Z)=H^{1,1}(V)\cap H^2(V;\C)$. Hence $H^2(V,\SO_V)=0$. 

Now, by (1) above, we know that
\[
\Pic(V) \cong H^2(V;\Z) \cong H^4(D;\Z).
\]
$\Pic(X)$ is isomorphic to the quotient group of $\Pic(V)$ modulo the subgroup generated by the classes of 
irreducible components of $D$. Hence $\Pic(X)=0$, and $X$ is factorial. Since there are no linear equivalence 
relations among the irreducible components of $D$, it follows that $A^*=\C^*$.
\QED

The same argument with the $\Z$-coefficients replaced by $\Q$-coefficients in the proof of Lemma \ref{Lemma 1.1} 
shows that {\em for a $\Q$-homology $n$-fold $X$ it holds that 
\begin{enumerate}
\item[(1)]
$X$ is $\Q$-factorial, i.e., $\Pic(X)$ is a finite abelian group, and $A^*=\C^*$.
\item[(2)]
$H^1(V,\SO_V)=H^2(V,\SO_V)=0$.
\end{enumerate}}
\svskip

The following result is an important consequence of a result of Hamm \cite{Hamm} that an affine 
variety of dimension $n$ defined over $\C$ has the homotopy type of a CW complex of real dimension $n$.

\begin{lem}\label{Lemma 1.2}
Let $X$ be an affine variety of dimension $n$. Then $H_n(X;\Z)$ is a free abelian group of finite rank. Furthermore, 
$H^i(X;\Z)=0$ for $i > n$.
\end{lem}
\Proof
The vanishing of $H^i(X;\Z)$ follows from the first assertion and the vanishing of $H_i(X;\Z)$ 
for $i > n$. By the universal coefficient theorem, we have 
\[
H^i(X:\Z)=\Hom(H_i(X;\Z),\Z)\oplus \Ext^1(H_{i-1}(X;\Z),\Z)\ .
\]
If $i=n+1$, $H^{n+1}(X;\Z)=0$ because $H_n(X;\Z)$ is torsion free. For $i > n+1$, the result is clear because 
$H_i(X;\Z)=H_{i-1}(X;\Z)=0$.
\QED

The following result (see \cite[Lemma 1.5]{Nori} for the proof) is frequently used below.

\begin{lem}\label{Lemma 1.3}
Let $f : Y \to X$ be a dominant morphism of algebraic varieties such that the general fibers are irreducible. 
Then the natural homomorphism $\pi_1(Y) \to \pi_1(X)$ is surjective. 
\end{lem}

As an extension of the above argument, one can give a different proof of the following result in \cite[Theorem 1]{M2}.
In the original proof, one has to use, in a crucial step of the proof, a rather difficult result that a smooth, 
quasi-affine surface $X$ of log Kodaira dimension $-\infty$ is either affine-ruled or contains $\A^2/\Gamma$ 
as an open set so that $X-\A^2/\Gamma$ is a disjoint union of affine lines, where $\Gamma$ is a small finite 
subgroup of $\GL(2,\C)$ (cf. the proof of \cite{M2}).

\begin{thm}\label{Theorem 1.4}
Let the additive group scheme $G_a$ act nontrivially on the affine $3$-space $\A^3$. Then $\A^3\quot G_a
\cong \A^2$.
\end{thm}

For the proof, we need two lemmas.

\begin{lem}\label{Lemma 1.5}
Let $Y$ be a smooth contractible affine threefold with a nontrivial $G_a$-action and let $X=Y\quot G_a$. 
Let $V$ be a normal projective surface containing $X$ as an open set such that $V$ is smooth along $D:=
V-X$ and $D$ is a divisor of simple normal crossings. Let $\rho : V' \to V$ be a minimal resolution of 
singularities of $X$. Let $p_g(V')$ be the geometric genus of the surface $V'$ which is a birational 
invariant independent of the choice of $V'$. Then the following assertions hold.
\begin{enumerate}
\item[(1)]
Both $X$ and $X^\circ$ are simply connected, where $X^\circ$ is the smooth part of $X$. 
\item[(2)]
$H_1(X;\Z)=0$. Further, if $p_g(V')=0$ and $H^1(D;\Z)=0$, then $H_2(X;\Z)=0$ and hence $X$ is contractible.
\item[(3)]
$X$ is smooth under the additional assumptions in the assertion {\em (2)}.
\end{enumerate}
\end{lem}
\Proof
(1)\ Write $Y=\Spec B$ and $X=\Spec A$. Then $B$ is factorial by Lemma \ref{Lemma 1.1} and hence $A$ is factorial, 
too. So, the singular locus of $X$ is a finite set of quotient singular points \cite[Lemma 3.4]{GMM}. 
Hence the quotient morphism $q : Y \to X$ has no fiber components of dimension two and $Y-q^{-1}(X^\circ)$ 
has dimension $\le 1$. Hence $\pi_1(q^{-1}(X^\circ))=\pi_1(Y)=(1)$. By Lemma \ref{Lemma 1.3}, both $X$ and $X^\circ$ 
are simply connected.

(2)\ Let $X'$ be the inverse image $\rho^{-1}(X)$ and let $E$ be the exceptional locus of $\rho$ which is 
a divisor with simple normal crossings. Let $F=D+E$. Then we have an exact sequence of 
cohomology groups
\[
H^3(V',F;\Z) \to H^3(V';\Z) \to H^3(F:\Z)
\]
where $H^3(F;\Z)=0$ and $H^3(V',F;\Z)\cong H_1(X^\circ;\Z)$ by the Lefschetz duality and $H_1(X^\circ;\Z)=0$
because $\pi_1(X^\circ)=1$. Hence $H_1(V';\Z)\cong H^3(V';\Z)=0$. As in the proof of Lemma \ref{Lemma 1.1}, 
this implies that $H^1(V';\C)=0$ and $H^1(V',\SO_{V'})=0$. Hence we have an exact sequence
\[
0 \to H^1(V',\SO_{V'}^*) \to H^2(V';\Z) \to H^2(V',\SO_{V'})\ .
\]
Since $H^2(V',\SO_{V'})=0$ by the hypothesis, we have $\Pic(V')\cong H^1(V',\SO_{V'}^*)$ $\cong H^2(V';\Z)$.
Namely the group $H^2(V';\Z)$ of topological $2$-cocycles is generated by algebraic classes of divisors. 
Since $X$ is factorial, $\Pic(V')$ is generated by the classes of irreducible components of $F$. 
Since $H^2(F;\Z)$ is a free abelian group generated by the classes of irreducible components of $F$, 
the natural homomorphism $H^2(V';\Z) \to H^2(F;\Z)$ is an isomorphism. 

Consider the long exact sequence associated to the pair $(V',F)$
\[
H^1(V';\Z) \to H^1(F;\Z) \to H^2(V',F;\Z) \to H^2(V';\Z) \to H^2(F;\Z),
\]
where $H^1(V';\Z)=0$ as $H_1(V';\Z)=0$ and $H^2(V';\Z)\to H^2(F;\Z)$ is an isomorphism. Hence $H^1(F;\Z) 
\cong H_2(V',F;\Z)$. Since the divisors $D$ and $E$ are disjoint from each other, $H^1(F;\Z)\cong 
H^1(D;\Z)\oplus H^1(E;\Z)$, where $H^1(D;\Z)=0$ by the hypothesis and $H^1(E;\Z)=0$ because $X$ has at worst 
quotient singular points. So, $H^1(F;\Z)=0$. By the Lefschetz duality, $H^2(V,F;\Z)\cong H_2(X^\circ;\Z)$. 
It follows that $H_2(X^\circ;\Z)=0$. 

Now let $\{P_1,\ldots,P_r\}$ be the set of singular points of $X$. Choose a closed neighborhood $T_i$ of $P_i$ 
for every $i$ so that $T_i\cap T_j=\emptyset$ if $i \ne j$. Let $\partial T_i$ be the boundary of $T_i$ and let 
$\partial T=\partial T_1\cup \cdots \cup \partial T_r$. By \cite[Proof of Lemma 2.2]{MS}, we have an exact 
sequence
\[
0 \to H_2(X^\circ;\Z) \to H_2(X;\Z) \to H_1(\partial T;\Z) \to H_1(X^\circ;\Z),
\]
where $H_1(X^\circ;\Z)=0$ as $\pi_1(X^\circ)=1$ and $H_1(\partial T;\Z)$ is a finite abelian group because 
$\pi_1(\partial T_i)$ is a finite group as the local fundamental group $\pi_{1,P_i}(X)$ of the quotient singular 
point $P_i$. Since $H_2(X^\circ;\Z)=0$ as shown above, $H_2(X;\Z)$ is a finite abelian group. Since $H_2(X;\Z)$ 
is torsion free by Lemma \ref{Lemma 1.2}, it follows that $H_2(X;\Z)=0$. 

(3)\ Since $\pi_1(X)=\pi_1(X^\circ)=1$ and $X$ is contractible, $X$ is smooth by the so-called {\em affine Mumford 
theorem} \cite[Theorem 3.6]{GKMR} which we state below.
\QED

\begin{lem}\label{Lemma 1.6}
Let $X$ be a normal affine surface such that $X$ is contractible and the smooth part $X^\circ$ of $X$ is simply 
connected. Then $X$ is smooth.
\end{lem}

Now the different proof of Theorem \ref{Theorem 1.4} is given as follows. With the notations in Lemma \ref{Lemma 1.5},
let $Y=\A^3$. Choose a general linear hyperplane $L$ in such a way that the quotient morphism $q$ restricted on $L$ 
is a dominant morphism to $X:=\A^3\quot G_a$ and $L$ meets the inverse image $q^{-1}(\Sing X)$ in finitely many points. 
Then it follows that $X$ is rational and $\lkd(X^\circ)=-\infty$. This implies that $p_g(V')=0$ and $H^1(D;\Z)=0$ 
with the notations in Lemma \ref{Lemma 1.5}. Hence $X$ is smooth. It is clear that $A$ is factorial, $A^*=\C^*$ and 
$\lkd(X)=-\infty$. By a characterization of the affine plane, $X$ is isomorphic to $\A^2$.

\begin{remark}\label{Remark 1.7}{\em 
To be accurate, we have to use the result in \cite[p. 49]{M2} mentioned before Theorem \ref{Theorem 1.4} in the course of 
the proof of Lemma \ref{Lemma 1.6}. Hence it is very desirable to give a new proof of Lemma \ref{Lemma 1.6} which is 
more topological and not depending on the result used in \cite{M2}. 

In \cite[Theorem 2.7]{KS}, Kaliman and Saveliev state a result stronger than Lemma \ref{Lemma 1.5}. They state that 
{\em for every nontrivial $G_a$-action on a smooth contractible affine algebraic threefold $Y$, the quotient 
$X=Y\quot G_a$ is a smooth contractible affine surface.}}
\end{remark}

\section{$\A^1_*$-fibration}
Let $p : Y \to X$ be a dominant morphism of algebraic varieties. We call $p$ an $\A^1$-(resp. $\A^1_*$-) fibration 
if general fibers of $p$ are isomorphic to $\A^1$ (resp. $\A^1_*$). Here $\A^1_*$ is the affine line with one point 
deleted. A {\em singular} fiber of $p$ is a fiber which is not scheme-theoretically 
isomorphic to $\A^1$ (resp. $\A^1_*$). We say that an $\A^1_*$-fibration is {\em untwisted} (resp {\em twisted})
if the generic fiber $Y_K:=Y\times_X\Spec K$ has two $K$-rational places at infinity (resp. if $Y_K$ has one non 
$K$-rational place at infinity), where $K$ is the function field of $X$ over $\C$. In the untwisted case, $Y_K$ 
is isomorphic to $\A^1_{*,K}:=\Spec K[t,t^{-1}]$, while $Y_K$ is a non-trivial $K$-form of $\A^1_{*,K}$ and 
$\Pic(Y_K) \cong \Z/2\Z$ in the twisted case. If $Y$ is a factorial affine variety, then any $\A^1_*$-fibration on 
$Y$ is untwisted because $Y_K$ is then factorial. In fact, if $Y=\Spec B$ and $X=\Spec A$ then $B$ is factorial 
and hence the quotient ring $B\otimes_AK$ is factorial. Since $Y_K=\Spec B\otimes_AK$, $Y_K$ is factorial. In the 
present article, an $\A^1_*$-fibration is {\em always assumed to be untwisted}.

If the additive group scheme $G_a$ (resp. the multiplicative group scheme $G_m$) acts on an affine threefold $Y$ 
then the algebraic quotient $Y \quot G_a$ (resp. $Y \quot G_m$) exists. Namely, if $Y=\Spec B$ and $X$ is the quotient 
of $Y$ by $G_a$ (resp. $G_m$), then the invariant subring $A=B^{G_a}$ (resp. $A=B^{G_m}$) is the coordinate ring 
of $X$ and the quotient morphism $q : Y \to X$ is given by the inclusion $A \hookrightarrow B$. Note that $A$ is 
finitely generated over $\C$ by a lemma of Zariski in the case of $G_a$ \cite{Z} and by the well-known result on 
reductive group actions in the case of $G_m$. In the case of $G_m$-actions, we assume {\em unless otherwise mentioned} 
that $\dim A =\dim B-1$. Hence $q$ has relative dimension one. In the case of a $G_a$-action, the quotient morphism 
$q : Y \to X$ is an $\A^1$-fibration which is not necessarily surjective, but in the case of a $G_m$-action, $q$ is 
surjective and either an $\A^1$-fibration or an $\A^1_*$-fibration according as general orbits admit fixed points 
or not. If the quotient morphism $q : Y \to X$ is an $\A^1_*$-fibration, it is untwisted by a theorem of Rosenlicht 
\cite{Rosenlicht}. We recall the following fundamental result on the quotients under reductive algebraic groups \cite{Borel}. 
When we consider the quotient morphism by $G_a$ (or $G_m$), we denote it by $q$ in most cases, but an 
$\A^1$-fibration (or $\A^1_*$-fibration) by $p$. 

\begin{lem}\label{Lemma 2.1}
Let $G$ be a connected reductive algebraic group acting on a smooth affine algebraic variety $Y$ and let $q : Y \to
X$ be the quotient morphism. Then, for two points $P_1, P_2$, we have $q(P_1)=q(P_2)$ if and only if 
$\ol{GP_1}\cap \ol{GP_2} \ne \emptyset$. Hence for any point $Q$ of $X$, the fiber $q^{-1}(Q)$ contains a unique 
closed orbit.
\end{lem}

In order to distinguish the case of surfaces from the case of threefolds, we use the notation like $p : X \to C$ 
in the case of surfaces and retain the notation $p : Y \to X$ for the case of threefolds or varieties in general.  

\begin{lem}\label{Lemma 2.2}
Let $X$ be a normal affine surface and let $p : X \to C$ be an $\A^1_*$-fibration. Then the following assertions 
hold.
\begin{enumerate}
\item[(1)]
A singular fiber $F$ of $p$ is written in the form $F=\Gamma+\Delta$, where $\Gamma=m\A^1_*, \emptyset$ or 
$m\A^1+n\A^1$ with two $\A^1$'s meeting in one point which may be a cyclic quotient singular point of $X$, 
$\Gamma \cap \Delta=\emptyset$ and $\Delta$ is a disjoint union of affine lines with multiplicities, each of 
which may have a unique cyclic singular quotient point of $X$. Each type of a singular fiber is realizable. 
\item[(2)]
If $p$ is the quotient morphism of a $G_m$-action on $X$ and $p$ is an $\A^1_*$-fibration, the part $\Delta$ is 
absent in a singular fiber $F$. Hence $F$ is either $m\A^1_*$ or $m\A^1+n\A^1$ in the above list.
\end{enumerate}
\end{lem}
\Proof
(1)\ In \cite[Lemma 4]{Mi8} and \cite[Proposition 5.1 (b)]{Zaidenberg}, the case where $X$ is a smooth affine surface 
is treated. In \cite[Lemma 2.9]{MS}, singular fibers of an $\A^1_*$-fibration are classified, where singularities of 
the surface $X$ are assumed to be quotient singularities, but one can drop this condition and make the same 
argument by assuming simply that $X$ is normal. The normal case is also treated in \cite[Proposition 3.8 and 
Theorem 4.18]{FZ}.

(2)\ By the assertion (1) and Lemma \ref{Lemma 2.1}, the only possible cases of a singular fiber of the quotient 
morphism $p$ are the two cases listed in the statement and the case where $F_\red$ is irreducible and isomorphic to 
$\A^1$ which may contain a unique cyclic quotient singularity. We show that the last case does not occur. Suppose 
that $F_\red \cong \A^1$. Then $F$ contains a unique fixed point $P$. Suppose first that $X$ is smooth at $P$. Look 
at the induced tangential representation of $G_m$ on $T_{X,P}$ which is diagonalizable and has weights $-a, b$
with $a, b > 0$. Then it follows that $F$ is locally near $P$ a union of two irreducible components meeting in $P$. 
This is a contradiction. Suppose that $X$ is singular at $P$. Let $\sigma : \wt{X} \to X$ be a minimal resolution 
of singularity at $P$. Then the inverse image $\sigma^{-1}(F_\red)$ is the proper transform $G$ of $F_\red$ and 
a linear chain of $\BP^1$'s with $G$ meeting one of the terminal components of the linear chain. Since the 
$G_m$-ation on $X$ lifts to $\wt{X}$, each irreducible component of $\sigma^{-1}(F_\red)$ is $G_m$-stable and 
the other terminal component, say $\wt{F}$, of the linear chain has an isolated fixed point, say $\wt{P}$. Note that 
we obtain an affine surface from $\wt{X}$ by removing all components of $\sigma^{-1}(F_\red)$ except for the 
component $\wt{F}$. Now looking at the induced tangential representation of $G_m$ on $T_{\wt{X},\wt{P}}$, we have 
a contradiction as in the above smooth case. So, the case $F_\red \cong \A^1$ does not occur.
\QED 

Special attention has to be paid in the case where the $G_m$-quotient morphism $q : Y \to X$ is an $\A^1$-fibration.
We consider the surface case first and then treat the case of threefolds.

\begin{lem}\label{Lemma 2.3}
Let $X$ be a normal affine surface with a $G_m$-action. Suppose that the quotient morphism $\rho : X \to C$ is an 
$\A^1$-fibration, where $C$ is a normal affine curve. Then the following assertions hold.
\begin{enumerate}
\item[(1)]
There exists a closed curve $\Gamma$ on $X$ such that the restriction of $q$ onto $\Gamma$ induces an isomorphism 
between $\Gamma$ and $C$, i.e., $\Gamma$ is a cross-section of $\rho$.
\item[(2)]
Suppose further that $X$ is smooth. Then every fiber of $q$ is reduced and isomorphic to $\A^1$. Hence $X$ is 
a line bundle over $C$ .
\item[(3)]
Suppose that $X$ is smooth and $C$ is rational. Then $X$ is isomorphic to a direct product $C\times \A^1$.
\end{enumerate}
\end{lem}
\Proof
(1)\ Let $\ol{\rho} : \ol{X} \to \ol{C}$ be a $G_m$-equivariant completion of $\rho : X \to C$. Namely, $\ol{X}$ 
is a normal projective surface containing $X$ as an open set, $\ol{C}$ is the smooth completion of $C$. We may 
assume that $\ol{X}$ is smooth along $\ol{X}\setminus X$, $\rho$ extends to a morphism $\ol{\rho}$ and the 
$G_m$-action extends to $\ol{X}$ so that $\ol{\rho}$ is a $\BP^1$-fibration. Since $\rho$ is an $\A^1$-fibration, 
there exists an open set $U$ of $C$ such that 
$\rho^{-1}(U) \cong U\times \A^1$ with $G_m$ acting on $\A^1$ in a natural fashion. Let $\Gamma_0$ be the fixed point 
locus in $\rho^{-1}(U)$ and let $\Gamma$ (resp. $\ol{\Gamma}$) be the closure of $\Gamma_0$ in $X$ (resp. $\ol{X}$). 
Then the restriction of $\ol{\rho}$ onto $\ol{\Gamma}$ is a birational morphism onto $\ol{C}$, hence an isomorphism. 
Thus $\ol{\Gamma}$ is a cross-section of $\ol{\rho}$. On the other hand, $\ol{\rho}$ has another cross-section 
$\ol{\Gamma}_{\infty}$ lying outside $X$. Since $\rho$ is an $\A^1$-fibration, every fiber of $\rho$ is a disjoint 
union of the affine lines, hence it consists of a single affine line by Lemma \ref{Lemma 2.1}. Suppose that 
for a point $\alpha \in C$, the fiber $F_\alpha$ does not meet $\Gamma$. Then $\ol{\Gamma}$ meets the fiber 
$\ol{F}_\alpha:=\ol{\rho}^{-1}(\alpha)$ outside $F_\alpha$, and $F_\alpha$ has an isolated $G_m$-fixed point $P$. 
This leads to a contradiction if we argue as in the proof of Lemma \ref{Lemma 2.2}. Namely, if $X$ is smooth at $P$, 
then consider the induced tangential representation at $P$. Otherwise, consider a minimal resolution of singularity 
at $P$ and a $G_m$-fixed point appearing in one of the terminal components of the exceptional locus which does not 
meet the proper transform of $F_\alpha$. Thus $\Gamma$ meets every fiber of $\rho$ and is smooth. 

(2)\ With the above notations, suppose that $X$ has a singular point on the fiber $F_\alpha$. Since $G_m$ acts 
transitively on $F_\alpha-F_\alpha\cap\Gamma$, $X$ has cyclic singularity at the point $P_\alpha:=F_\alpha\cap \Gamma$. 
Let $\sigma : \wt{X} \to \ol{X}$ be the minimal resolution of singularities and let $\Delta=\sigma^{-1}(P_\alpha)$. 
The $\BP^1$-fibration $\ol{\rho} : \ol{X} \to \ol{C}$ lifts to a $\BP^1$-fibration $\wt{\rho} : \wt{X} \to \ol{C}$. 
The proper transform $\wt{\Gamma}$ of $\ol{\Gamma}$ is a cross-section of $\wt{\rho}$. The fiber $\wt{F}_\alpha$ 
of $\wt{\rho}$ corresponding to $\ol{F}_\alpha$ contains the linear chain $\Delta$ in such a way that one terminal 
component, say $G$, meets $\wt{\Gamma}$ and the other terminal component meets the proper transform of $F_\alpha$. Further, 
all the components except for the terminal component $G$ meeting $\wt{\Gamma}$ can be contracted smoothly. In fact, 
we can replace the fiber $F_\alpha$ by $G-\{Q\}$ without losing the affineness of $X$, where $Q$ is the point 
where $G$ meets the adjacent component of $\wt{F}_\alpha$. If $X$ is smooth along $F_\alpha$, it is reduced. If 
$X$ is smooth, the $\A^1$-fibration $\rho$ has no singular fibers. Hence $X$ is an $\A^1$-bundle over $C$ 
\cite[Lemma 1.15]{GMM}. Since the two cross-sections $\ol{\Gamma}$ and $\ol{\Gamma}_\infty$ are disjoint from each 
other over $C$, it follows that $X$ is a line bundle over $C$.

(3)\ Since $\Pic(C)=0$, every line bundle is trivial. 
\QED     
 
The proof of the assertion (2) implies that every normal affine surface with a $G_m$-action and an $\A^1$-fibration 
as the quotient morphism is constructed from a line bundle over $C$ by a succession of blowing-ups with centers on the 
zero section and contractions of linear chains. The following result is well-known.

\begin{lem}\label{Lemma 2.4}
Let $Y=\Spec B$ be a smooth affine threefold with a $G_m$-action and let $q : Y \to X:=\Spec A$ be the quotient 
morphism. Assume that $q$ is an $\A^1$-fibration and that $q$ is equi-dimensional. Then the following assertions hold.
\begin{enumerate}
\item[(1)]
We may assume that $B$ is positively graded, that is, $B$ is a graded $A$-algebra indexed by $\Z_{\ge 0}$.
\item[(2)]
$A$ is normal and factorially closed in $B$.
\item[(3)]
$q : Y \to X$ has a cross-section $S$. Namely, each fiber of $q$ has dimension one and meets $S$ in one point 
transversally.
\item[(4)]
$X$ is smooth, and $Y$ is a line bundle over $X$.
\item[(5)]
If $Y$ is factorial, then the equi-dimensionality condition is automatically satisfied and $Y$ is isomorphic to 
a direct product $X\times \A^1$ with $G_m$ acting on the factor $\A^1$ in a natural fashion.
\end{enumerate}
\end{lem}
\Proof
(1)\ By \cite[Lemma 1.2]{GMM}, there exists a nonzero element $a \in A$ such that $B[a^{-1}]=A[a^{-1}][u]$, where $u$ is an 
element of $B$ algebraically independent over the quotient field $Q(A)$. In fact, $u$ is a variable on a general 
fiber of $q$. We may assume that the $G_m$-action on the fiber is given by ${}^tu=tu$. Note that the $G_m$-action 
gives a $\Z$-graded ring structure on $B$. Let $b$ be a homogeneous element of $B$. We can write $a^rb=f(u)$, 
where $f(u)=a_0u^m+a_1u^{m-1}+\cdots+a_m \in A[u]$ with $a_0 \ne 0$. Then we have
\[
a^r({}^tb)=f(tu)=a_0t^mu^m+a_1t^{m-1}u^{m-1}+\cdots+a_m\ .
\]
Hence $a^rb=a_0u^m$ and $b$ has degree $m$. This implies that every homogeneous element of $B$ has degree $\ge 0$.

(2)\ The normality of $A$ is well-known \cite[p.100]{Kraft}, and the factorial closedness follows from the assertion (1). 

(3)\ Let $S_0$ be the closed set in the open set $q^{-1}(D(a))$ defined by $u=0$. Hence $S_0$ is the locus of the 
fixed points in the fibers $q^{-1}(Q)$ when $Q \in D(a)$. Let $S$ be the closure of $S_0$ in $Y$. Then $S$ is a 
rational cross-section of $q$. We shall show that $S$ meets every fiber of $q$ in one point. Let $P$ be a point of $X$,  
let $C$ be a general irreducible curve on $X$ passing through $P$ and let $Z=Y \times_X C$. Then $q_C : Z \to C$, the 
projection onto $C$, is an $\A^1$-fibration. Let $\nu_Z : \wt{Z} \to Z$ and $\nu_C : \wt{C} \to C$ be the normalizations 
of $Z$ and $C$ respectively. Then there exists an $\A^1$-fibration $\rho : \wt{Z} \to \wt{C}$ such that 
$\nu_C\cdot\rho = q_C\cdot \nu_Z$. The normal affine surface $\wt{Z}$ has the induced $G_m$-action and $\rho$ is the 
quotient morphism. Choosing the open set $D(a)$ small enough, we may assume that $C \cap D(a)$ is a non-empty set 
contained in the smooth part of $C$. Let $\Gamma_0=S_0\cap Z$ and let $\wt{\Gamma}$ be the closure of $\Gamma_0$ in 
$\wt{Z}$. Let $\wt{P}$ be a point of $\wt{C}$ lying over $P$. Let $\wt{F}=\rho^{-1}(\wt{P})$ and $F=q^{-1}(P)$. 
By Lemma \ref{Lemma 2.3}, $\wt{\Gamma}$ meets $\wt{F}$ in a single point. This implies that $S$ meets $F$ in a point,  
say $R$. By \cite[Theorem 1.15 and Lemma 3.5]{GMM} and Lemma \ref{Lemma 2.1}, $F$ consists of one irreducible component 
which is contractible and smooth outside the point $R$ because it admits a non-trivial $G_m$-action with a fixed point $R$. 
Since $Y$ is smooth, the local intersection multiplicity $i(S,F;R)=1$. This implies that $F$ is reduced and smooth also 
at $R$. Hence $F$ is isomorphic to $\A^1$. 

(4)\ The morphism $q : Y \to X$ induces a quasi-finite birational morphism $q|_S : S \to X$. Then it is an isomorphism 
by Zariski's main theorem. Since $S$ is smooth by the argument in the proof of the assertion (3), $X$ is also smooth. 
Furthermore, every fiber of $q$ is an affine line. Hence $Y$ is an $\A^1$-bundle over $X$. Since a smooth completion 
$\ol{q} : \ol{Y} \to \ol{X}$ is a $\BP^1$-fibration with two cross-sections (the zero section and the infinity section) 
which are disjoint over $X$, $Y$ is a line bundle over $X$.

(5)\ By \cite[Lemma 1.10]{GMM}, $q$ does not have codimension one fiber components. Since $X$ is factorial, 
$Y$ is isomorphic to $X\times \A^1$. This fact follows from the proof of the assertion (1). In fact, one can take the 
element $u$ to be a prime element of $B$. Then, the relation $a^rb=a_0u^m$ for a homogeneous element $b$ of degree $m$ 
implies that $b=cu^m$ for a certain element $c \in A$. Hence $B=A[u]$. 
\QED

In \cite[Theorem 1.4]{GMM}, it is shown that any $\A^1$-fibration $p : Y \to X$ on an affine threefold $Y$ is factored by 
the quotient morphism $q : Y \to Z$ by a certain $G_a$-action on $Y$. Meanwhile, this is not the case with 
an $\A^1_*$-fibration. Examples can be easily produced by using Lemma \ref{Lemma 2.2}, (1). 

A homology threefold $Y$ is a {\em contractible threefold} if it is topologically contractible. A homology threefold is 
contractible if and only if $\pi_1(Y)=1$. We consider a $G_m$-action on such a threefold.

\begin{lem}\label{Lemma 2.5}
Let $G_m$ act nontrivially on a $\Q$-homology threefold $Y$. Then the fixed point locus $Y^{G_m}$ is a non-empty connected  
closed subset of $X$. Furthermore, $H_i(Y^{G_m};\Z/p\Z)=0$ for every $i >0$ and for infinitely many prime numbers $p$. 
\end{lem}
\Proof
Since $H_i(Y;\Z)$ is a finite group, there exist infinitely many primes $p$ such that $Y$ is $\Z/p\Z$-acyclic, i.e., 
$H_i(Y;\Z/p\Z)=0$ for every $i > 0$. We choose such a prime $p$ with the additional property that $p$ does not divide 
any weight of the induced tangential $G_m$-action on $T_{Y,Q}$ for every fixed point $Q$. We take a point $Q \in Y^{G_m}$.
By \cite[Theorem 1]{Koras}, there exists a $G_m$-stable open neighborhood $U$ of $Q$ in the Euclidean topology and algebraic 
functions $x,y,z$ around the point $Q$ which form a system of analytic coordinates on $U$. The $G_m$-action is given by 
${}^t(x,y,z)=(t^ax,t^by,t^cz)$ with integers $a, b, c$, where $p \nmid \alpha$ for any $\alpha \in \{a, b,c\}\setminus\{0\}$. 
Then $U\cap Y^{G_m}=U\cap Y^{\Z/p\Z}$. Since $Y^{G_m} \subseteq Y^{\Z/p\Z}$, it follows that $Y^{G_m}$ is a connected 
component of $Y^{\Z/p\Z}$. By \cite{Floyd}, we have the inequality
\[
\sum_i \dim H_i(Y^{\Z/p\Z};\Z/p\Z) \le \sum_i \dim H_i(Y; \Z/p\Z) =1.
\]
This implies that $Y^{\Z/p\Z}$ is connected and equal to $Y^{G_m}$. Hence $Y^{G_m}$ is connected and $H_i(Y^{G_m};\Z/p\Z)=0$ 
for every $i >0$. 
\QED

We have the following result.

\begin{thm}\label{Theorem 2.6}
Let $G_m$ act nontrivially on a $\Q$-homology threefold $Y$ and let $q : Y \to X$ be the quotient morphism. Suppose 
that $q$ has relative dimension one. Then the following assertions hold.
\begin{enumerate}
\item[(1)]
$Y^{G_m}$ is $\Q$-acyclic. If $Y$ is a homology threefold, then $Y^{G_m}$ is $\Z$-acyclic.
\item[(2)]
If $\dim Y^{G_m}=2$ then $Y^{G_m} \cong X$ and $Y$ is a line bundle over $X$. If $Y$ is a homology threefold, then 
$Y \cong X\times \A^1$. 
\item[(3)]
If $\dim Y^{G_m}=1$ then $Y^{G_m} \cong \A^1$ and $q|_{Y^{G_m}} : Y^{G_m} \to X$ is a closed embedding. If $Y$ is 
contractible and $q$ is equi-dimensional, then $X$ is a smooth contractible surface of log Kodaira dimension $-\infty$ 
or $1$. 
\item[(4)]
If $\dim Y^{G_m}=0$ then the $G_m$-action on $Y$ is hyperbolic, i.e., the tangential representation on $T_{Y,Q}$ for the 
unique fixed point $Q$ has mixed weights, e.g., either $a_1 < 0, a_2 > 0$ and $a_3 > 0$, or $a_1 >0, a_2 < 0$ 
and $a_3 < 0$. In this case, there is an irreducible component of codimension one contained in the fiber of $q$ 
passing through the point $Q$.
\end{enumerate}
\end{thm}
\Proof
(1)\ The $\Q$-acyclicity is verified below case by case according to $\dim Y^{G_m}$. When $X$ is a homology threefold, 
the acyclicity of $Y^{G_m}$ follows from the Smith theory \cite[\S 22.3]{R1}. 

(2)\ By Lemma \ref{Lemma 2.1}, $Y^{G_m}$ lies horizontally to the morphism $q$ and the restriction $q|_{Y^{G_m}} : 
Y^{G_m} \to X$ is a bijection (hence a birational morphism). In particular, $Y^{G_m}$ is irreducible. Further, 
$q : Y \to X$ is an $\A^1$-fibration. Since $Y$ is $\Q$-factorial by a remark after Lemma \ref{Lemma 1.1}, $q$ is 
equi-dimensional (cf. \cite[Lemma 1.14]{GMM}). Then Lemma \ref{Lemma 2.4} implies that $X$ is smooth and $q : Y \to X$ 
is a line bundle. By Zariski's main theorem, we conclude that $Y^{G_m} \cong X$. Furthermore, since $Y$ is contractible 
to $X$, it follows that $Y^{G_m}$ is $\Q$-acyclic (resp. $\Z$-acyclic) if $Y$ is a $\Q$-homology threefold 
(resp. $\Z$-homology threefold). If $Y$ is a homology threefold, then $X$ is factorial and a line bundle over $X$ is 
trivial. Hence $Y \cong X \times \A^1$. 

(3)\ $Y^{G_m}$ is a connected curve. Further, $Y^{G_m}$ is smooth. In fact, the smoothness of the fixed point locus is 
a well-known fact for a reductive algebraic group action on a smooth affine variety (cf. \cite{R1})\footnote
{In the case where $G$ is a connected reductive algebraic group acting on a smooth affine variety $Y$, the $G$-action 
near a fixed point $Q$ is locally analytically $G$-equivalent to a linear representation \cite{Koras, Luna}. In the case of 
a linear representation of $G$ on the affine space $\C^n$, the fixed point locus near the origin is a linear subspace 
and it is $G$-equivariantly a direct summand of $\C^n$ by the complete reducibility of $G$. The smoothness of the fixed 
point locus $Y^G$ near the point $Q$ follows from this observation. }. Hence $Y^{G_m}$ is irreducible. Then $Y^{G_m}$ is 
an affine line because $H_1(Y^{G_m};\Z/p\Z)=0$. By Lemma \ref{Lemma 2.1}, $q$ induces a closed embedding of $Y^{G_m}$ into 
$X$ \footnote{Let $Q \in Y^{G_m}$. In view of \cite{Koras}, there exists a system of local (analytic) coordinates 
$\{x,y,z\}$ at $Q$ such that $G_m$ acts as ${}^t(x,y,z)=(x,t^{-a}y,t^bz)$ with $ab > 0$ and $y=z=0$ defining the curve 
$Y^{G_m}$ near $Q$. Then the quoteint surface $X$ at the point $P=q(Q)$ has a system of local analytic coordinates 
$\{x, x^{b'}y^{a'}\}$ with $a'=a/d, b'=b/d$ and $d=\gcd(a,b)$. Then the curve $q(Y^{G_m})$ is defined by $x^{a'}y^{b'}=0$ 
near $P$ and $X$ is smooth near the curve $q(Y^{G_m})$. Hence $Y^{G_m}$ and $q(Y^{G_m})$ are locally isomorphic at $Q$ 
and $P$. Since $q|_{Y^{G_m}} : Y^{G_m} \to X$ is injective, it is a closed embedding. Furthermore, the fiber of $q$ 
through the point $Q$ is a cross $a'\A^1+b'\A^1$. See the definition after Corollary \ref{Corollary 2.8}.}. 
If $Y$ is contractible, $X$ is contractible by \cite[Theorem B]{KPR}. By Lemma \ref{Lemma 1.3}, $\pi_1(X^\circ)=(1)$, where 
$X^\circ$ is the smooth part of $X$. Then $X$ is smooth by Lemma \ref{Lemma 1.6}. So, $X$ is a smooth contractible surface 
containing a curve isomorphic to $\A^1$. Hence $X$ has log Kodaira dimension $-\infty$ or $1$ (cf. \cite{Zaidenberg, GM}).

(4)\ The induced $G_m$-action on the tangent space $T_{Y,Q}$ of the unique fixed point $Q$ must have mixed weights. 
Otherwise, $\dim Y^{G_m} > 0$ near the point $Q$. This is a contradiction. To prove the last assertion, we have only to 
show it when $Y$ is idetified with the tangent linear $3$-space $T_{Y,Q}$. Then the $G_m$-action is diagonalized and 
${}^t(x,y,z)=(t^{a_1}x,t^{a_2}y,t^{a_3}z)$ with respect to a suitable system of coordinates $(x,y,z)$. Then the invariant 
elements, viewed as elements in $\Gamma(Y,\SO_Y)$, are divisible by $x$. Hence the fiber $F$ of $q$ passing through $Q$ 
contains an irreducible component $\{x=0\}$ which is a hypersurface of $Y$. 
\QED

\begin{question}\label{Question 2.7}{\em
In the case (3) above, is $Y$ isomorphic $G_m$-equivariantly to a direct product $Z \times Y^{G_m}$, where $Z$ 
is a smooth affine surface with $G_m$ acting on it? If this is the case, $q$ is the direct product $q=q'\times Y^{G_m}$ 
of the quotient morphism $q' : Z \to C$, where $C=Z\quot G_m$.}
\end{question}

We discuss this question in Section $4$ in detail. Theorem \ref{Theorem 2.6} has the following consequence.

\begin{cor}\label{Corollary 2.8}
Let the notations and the assumptions be the same as in Theorem \ref{Theorem 2.6}. If the quotient morphism $q : Y \to X$ 
has equi-dimension one, then the fixed point locus $Y^{G_m}$ has positive dimension.
\end{cor}

Now we consider an $\A^1_*$-fibration $p : Y \to X$, where $Y$ is a smooth affine threefold and $X$ is a normal affine 
surface. For a point $P$, the scheme-theoretic fiber $Y\times_X\Spec k(P)$ is denoted by $p^{-1}(P)$ or $F_P$, where 
$k(P)$ is the residue field of $P$ in $X$. A fiber $F_P$ of $p$ is called {\em singular} if $F_P$ is not isomorphic to 
$\A^1_*$ over $k(P)$. The set of points $P \in X$ such that $F_P$ is singular is denoted by $\Sing^{st}(p)$ and called 
the {\em strict singular locus} or the {\em degeneracy locus} of $p$. For a technical reason (cf. Lemma \ref{Lemma 2.13} 
and a remark below it), we define the {\em singular locus} $\Sing(p)$ as the union of $\Sing^{st}(p)$ and the set 
$\Sing(X)$ of singular points. A singular fiber $F_P$ is called a {\em cross} (resp. {\em tube}) \footnote{The naming of 
cross is apparent. The fiber $m\A^1_*$ is a projective line with two end points lacking and thickened with multiplicity. 
It looks like a tube.} if it has the form $F_P \cong m\A^1+n\A^1$ (resp. $m\A^1_*$). When we speak of a cross on a threefold 
or a surface which is smooth at the intersection point $Q$ of two lines, we assume that two affine lines meet each other 
transversally at $Q$. By abuse of the terminology, we call a singular fiber $F$ on a normal surface a cross even when 
the surface has a cyclic quotient singularity at the intersection point $Q$ and the proper transforms of the two lines 
on the minimal resolution of singularity form a linear chain together with the exceptional locus. 

We first recall a result of Bhatwadekar-Dutta \cite[Theorem 3.11]{BD}.

\begin{lem}\label{Lemma 2.9}
Let $R$ be a Noetherian normal domain and let $A$ be a finitely generated flat $R$-algebra such that 
\begin{enumerate}
\item[(1)]
The generic fiber $K\otimes_R A$ is a Laurent polynomial ring $K[T,T^{-1}]$ in one variable over $K=Q(R)$.
\item[(2)] 
For each prime ideal $\gp$ of $R$ of height one, the fiber $\Spec A\otimes_Rk(\gp)$ is geometrically integral but is not 
$\A^1$ over $k(\gp)$.
\end{enumerate}
Then there exists an invertible ideal $I$ in $R$ such that $A$ is a $\Z$-graded $R$-algebra isomorphic to the $R$-subalgebra 
$R[IT, I^{-1}T^{-1}]$ of $K[T,T^{-1}]$. In particular, $\Spec A$ is locally $\A^1_*$ and hence an $\A^1_*$-fibration over 
$\Spec R$.
\end{lem}

In geometric terms, a weaker version of Lemma \ref{Lemma 2.9} in the case of dimension three is stated as follows.

\begin{lem}\label{Lemma 2.10}
Let $p : Y \to X$ be an $\A^1_*$-fibration with a smooth affine threefold $Y$ and let $P$ be a closed point of $X$. Suppose 
that the following conditions are satisfied.
\begin{enumerate}
\item[(1)]
There is an open neighborhood $U$ of $P$ such that $U$ is smooth and $p$ is equi-dimensional over $U$. 
\item[(2)]
Every fiber $F_{P'}$ of $p$ for $P' \in U\setminus\{P\}$ is isomorphic to $\A^1_*$.
\end{enumerate}
Then the fiber $F_P$ is isomorphic to $\A^1_*$ and $p : p^{-1}(U) \to U$ is an $\A^1_*$-bundle. 
\end{lem}
\Proof
By replacing $X$ by a smaller affine open neighborhood of $P$ contained in $U$, we may assume that $X$ is smooth, $p$ 
is equi-dimensional and $\Sing(p)\cap (X\setminus\{P\})=\emptyset$. Then $p$ is flat. Further, there exists a $\BP^1$-fibration 
$\ol{p} : \ol{Y} \to X$ such that $Y$ is an open set of $\ol{Y}$ and $\ol{p}|_{Y}=p$. Since the $\A^1_*$-fibration $p$ 
is assumed to be untwisted, the generic fiber $\ol{Y}_K$ has two $K$-rational points $\xi_1, \xi_2$, where $K$ is 
the function field of $X$ over $\C$. Let $S_i$ be the closure of $\xi_i$ in $\ol{Y}$ for $i=1, 2$. By the assumption 
(2) above, $S_1, S_2$ are cross-sections over $X\setminus\{P\}$. Namely, $p^{-1}(X\setminus\{P\})$ is an $\A^1_*$-bundle. 
Hence $p^{-1}(X\setminus\{P\})\cong \CSpec (\SO_{X\setminus\{P\}}[\SL',{\SL'}^{-1}])$, where $\SL'$ is an invertible sheaf on 
$X\setminus\{P\}$. Since $X$ is smooth and $P$ is a point of codimension two, $\SL'$ is extended to an invertible sheaf 
$\SL$ on $X$. Since $Y$ is affine and $\codim_Yp^{-1}(P)=2$, it follows that $Y\cong \CSpec (\SO_X[\SL,\SL^{-1}])$. 
So, $p : Y \to X$ is an $\A^1_*$-bundle and $F_P \cong \A^1_*$.
\QED

Let $F_P$ be a fiber of an $\A^1_*$-fibration $p : Y \to X$ with $P \in X$. Let $C$ be a general smooth curve on $X$ passing 
through $P$. Let $Z$ be the normalization of $Y\times_XC$. Then $p$ induces an $\A^1_*$-fibration $p_C : Z \to C$ and the 
fiber $\wt{F}_P$ over the point $P$ is a finite covering of the fiber $F_P$ (in the sense that the normalization morphism 
induces a finite morphism $\wt{F}_P \to F_P$. In view of Lemma \ref{Lemma 2.2}, we say that the $\A^1_*$-fibration $p$ is 
{\em unmixed} if $\wt{F}_P$ for every $P \in X$ is $\A^1_*$, $\A^1+\A^1$, or a disjoint union of $\A^1$'s when taken with 
reduced structure. This definiton does not depend on the choice of the curve $C$. 

\begin{lem}\label{Lemma 2.11}
Let $p : Y \to X$ be an unmixed $\A^1_*$-fibration such that $Y$ and $X$ are smooth. Then the singular locus $\Sing(p)$ is a 
closed subset of $X$ of codimension one.
\end{lem}
\Proof
Lemma \ref{Lemma 2.9} or Lemma \ref{Lemma 2.10} implies that every irreducible component of the closure $\ol{\Sing(p)}$ of 
$\Sing(p)$ has codimension one. 

So, we have only to show that $\Sing(p)$ is a closed set. It suffices to show that if the fiber $F_P$ over $P \in X$ 
is $\A^1_*$ then there exists an open neighborhood $U$ of $P$ such that $F_{P'} \cong \A^1_*$ for every $P' \in U$. 
Let $C$ be a curve passing through $P$. If $C$ is not a component of $\ol{\Sing(p)}$, then the fiber of $p$ over 
the generic point of $C$ is $\A^1_*$ and hence geometrically integral. Suppose that $C$ is a component of $\ol{\Sing(p)}$. 
The unmixedness condition on $p$ implies that for a general point $P'$ of $C$, the fiber $F_{P'}$ is either irreducible 
and dominated by $\A^1_*$ (case (i)), or each irreducible component of $F_{P'}$ is dominated by a contractible curve 
(case (ii)). Suppose that the case (ii) occurs. Then there exists a normal affine surface $\wt{Z}$ with an $\A^1$-fibration 
$\wt{p} : \wt{Z} \to \wt{C}$ satisfying the following commutative diagram 
\[
\CD
\wt{Z}~~ @>\mu >> Y\times_XC \\
@VV\wt{p} V                  @VV p_C V \\
\wt{C}~~  @>\nu >> C
\endCD
\]
where $\mu$ is a quasi-finite morphism and $\nu$ is the normalization morphism. Since any singular fiber of an $\A^1$-fibration 
on a normal affine surface is a disjoint union of the affine lines, it turns out that the fiber $F_P$, which is a fiber of 
$p_C$ in the above diagram, is dominated by the affine line, which is a connected component of the fiber of $\wt{p}$. This 
is a contradiction, and the case (ii) cannot occur. Now removing from $X$ all irreducible components of $\ol{\Sing(p)}$ 
for whose general points the case (ii) occurs and replacing $Y$ by the inverse image of the open set of $X$ thus obtained, 
we may assume that the case (i) occurs for general points of the irreducible components of $\ol{\Sing(p)}$. Let 
\[
dp : \ST_Y \to p^*\ST_X
\]
be the tangential homomorphism of the tangent bundles on $Y$ and $X$. Let $\SC$ be the cokernel of $dp$. Then $\SC$ is a 
coherent $\SO_Y$-Module. Hence $\Supp(\SC)$ is a closed set $T$ such that $T=p^{-1}(p(T))$ and $T=\ol{\Sing(p)}$. Hence 
the point $P$ belongs to $p(T)$. However, $dp$ is everywhere surjective on the fiber $F_P$ and $F_P \not\subset \Supp(\SC)$. 
This is a contradiction. 

By the above argument, every irreducible curve $C$ through the point $P$ has the property that the general fibers are isomorphic 
to $\A^1_*$. Hence the fiber of $p$ over the generic point of $C$ is geometrically integral. By Lemma \ref{Lemma 2.9}, we know 
that $B\otimes_A\SO_P\cong \SO_P[x,x^{-1}]$. This implies that there exists an open neighborhood $U$ of $P$ satisfying 
the required property.
\QED

In the proofs of Lemmas \ref{Lemma 2.9}, \ref{Lemma 2.10} and \ref{Lemma 2.11}, the flatness condition on the morphism $p$ 
should not be overlooked. Hence the assumption that $X$ be normal instead of being smooth does not seem to be sufficient 
for the conclusion. So we ask the following question. 

\begin{question}\label{Question 2.12}{\em
Let $p : Y \to X$ be an $\A^1_*$-fibration with a smooth affine threefold $Y$. Is a point $P \in X$ smooth if 
the fiber $F_P$ is isomorphic to $\A^1_*$?}
\end{question}

The following is a partial answer to this question.

\begin{lem}\label{Lemma 2.13}
Let $Y$ be a smooth affine threefold with an effective $G_m$-action and let $q : Y \to X$ be the quotient morphism. Suppose 
that the $G_m$-action is fixed point free. Then $X$ has at worst cyclic quotient singular points.
\end{lem}
\Proof
Suppose that $P$ is a singular point of $X$. Let $Q$ be a point of $Y$ such that $q(Q)=P$ and let $G_Q$ be the isotropy group 
at $Q$ which is a finite cyclic group. Then, by Luna's \'etale slice theorem \cite{Luna}, there exists an affine subvariety 
$V$ of $Y$ with a $G_Q$-action and an \'etale morphism $\varphi/G_Q : V/G_Q \to X$ such that $Q \in V$, $V$ is smooth at $Q$ 
and $(\varphi/G_Q)(Q)=P$. This implies that $X$ has at worst cyclic quotient singularity at $P$. Since the ground field is 
$\C$, we can work with an analytic slice instead of an \'etale slice. 
\QED

We can easily show that if $X$ is smooth at $P$ then the fiber $F_P$ has multiplicity equal to $m=|G_Q|$ near the point $Q$. 
In fact, if $V$ is an analytic slice with coordinates $x,y$, the $G_Q$-action on $V$ is given by ${}^\zeta(x,y)=
(\zeta x,\zeta^by)$ with a generator $\zeta$ of $G_Q$ which is identified with an $m$-th primitive root of unity, where $0 \le 
b < m$. Then $P$ is singular if and only if $b > 0$. Hence $b=0$ if $P$ is smooth. This implies that $(x^m,y)$ is a local system 
of parameters of $X$ at $P$. Hence the fiber $F_P$ has multiplicity $m$. If $X$ is singular at $P$, the fiber $F_P$ is a multiple 
fiber. In fact, with the above notations, $\gm_{X,P}\wh{\SO}_{Y,Q}$ is generated by 
\[
\{x^m, y^m\}\cup\{x^sy^n \mid bn+s \equiv 0 \pmod{m},~0 <s < m,~0 < n < m\}.
\]
Since $\wh{\SO}_{Y,Q}=\C[[x,y,z]]$ with a fiber coordinate $z$ of $F_P$, it follows that $\gm_{X,P}\wh{\SO}_{Y,Q}$ is not 
a radical ideal. Hence the fiber $F_P$ is not reduced near the point $Q$. 
\svskip  

The singular locus $\Sing(p)$ may consist of a single point as shown in the following example.

\begin{example}\label{Example 2.14}{\em 
Let $\A^3 \to \A^2$ be the morphism defined by $(x,y,z) \mapsto (xy, x^2(xz+1))$, where $(x,y,z)$ is a system of 
coordinates of $\A^3$. Then $p$ is an $\A^1_*$-fibration and $\Sing(p)=\{(0,0)\}$. In fact, set $\alpha=xy,\ \beta=x^2(xz+1)$. 
Then $p^{-1}(\alpha,\beta)=\{(\alpha y^{-1},y,\alpha^{-3}y(\beta y^2-\alpha^2) \mid y \in \C^*\} \cong \A^1_*$ 
if $\alpha \ne 0$, $p^{-1}(0,\beta) = \{y=0,~x^2(xz+1)=\beta \} \cong \A^1_*$ if $\alpha=0,~\beta \ne 0$ and 
$p^{-1}(0,0) = \{x=0\}\cup\{y=xz+1=0\} \cong \A^2\cup \A^1_*$ if $\alpha=\beta=0$. \QED }
\end{example}

The following example shows that a disjoint union of affine lines may appear as a singular fiber of 
an $\A^1_*$-fibration. 

\begin{example}\label{Example 2.15}{\em
Let $p : \A^3 \to \A^2$ be the morphism defined by $(x,y,z) \mapsto (xy, x^2(yz+1))$. Then $p$ is an $\A^1_*$-fibration 
such that $\Sing(p) \cong \A^1$ and the fibers are given as follows. Set $\alpha=xy,\ \beta=x^2(yz+1)$. Then 
$p^{-1}(\alpha,\beta)=\{(\alpha y^{-1}, y, \alpha^{-2}(\beta y-\alpha^2 y^{-1})\mid y \in \C^*\}\cong \A^1_*$ if 
$\alpha \ne 0$, $p^{-1}(0,\beta)=\{y=0,~x^2=\beta\}\cong \A^1(x=\sqrt{\beta},y=0)\sqcup \A^1(x=-\sqrt{\beta},y=0)$ 
if $\alpha=0,~\beta \ne 0$ and $p^{-1}(0,0)=\{x=0\} \cong \A^2$ if $\alpha=\beta=0$. \QED}
\end{example}

Winkelmann \cite{W} kindly communicated us of the following example of a flat $\A^1_*$-fibration $p : \A^3 \to \A^2$ 
which is not surjective.  

\begin{prop}\label{Proposition 2.16}
Let $f(x)=\prod_{i=1}^n(x-i)$ and let $p : \A^3 \to \A^2$ be the morphism defined by $(x,y,z) \mapsto (x+zx+zyf(x), x+yf(x))$. 
Then the following assertions hold.
\begin{enumerate}
\item[(1)]
A general fiber is isomorphic to $\A^1_{n*}$, where $\A^1_{n*}$ is the affine line minus $n$ points. Hence $p$ is an 
$\A^1_*$-fibration if $n=1$.
\item[(2)]
$\im(p)=\A^2\setminus S$, where $S=\{(k,0)\mid 1 \le k \le n\}$.
\item[(3)]
$p$ is equi-dimensional. Hence $p$ is flat.
\item[(4)]
The singular locus $\Sing(P)$, that is the locus of points $P$ of the base $\A^2$ over which $p^{-1}(P) \not\cong \A^1_{n*}$, 
is the union of $\{(\alpha,0) \mid \alpha \in \C\}$ and $\bigcup_{k=1}^n\{(\alpha,k) \mid \alpha \in \C\}$. Hence $\Sing(p)$ 
consists of disjoint $n+1$ affine lines. For any $\alpha \in \C-\{1,\ldots,n\}$, the fiber $p^{-1}(\alpha,0) \cong \A^1$. 
For any $k \in \{1,\ldots,n\}$ and $\alpha \in \C$, the fiber $p^{-1}(\alpha,k)\cong \A^1\cup \A^1_{(n-1)*}$.
\end{enumerate}
\end{prop}
\Proof
(1)\ Set $\alpha=x+z(x+yf(x))$ and $\beta=x+yf(x)$. Then $\alpha=x+z\beta$. Hence if $\beta \ne 0$, we have 
$z=(\alpha-x)/\beta$. Further, if $\beta \not\in\{1,2,\ldots,n\}$, then the equation $\beta=x+yf(x)$ is solved as 
$y=(\beta-x)/f(x)$. In fact, if $(x,y)$ is a solution of $\beta=x+yf(x)$ and $f(x)=0$, then $x=\beta$ and $\beta 
\in \{1,2,\ldots,n\}$, a contradiction. Hence $p^{-1}(\alpha,\beta)\cong \A^1-\{1,2,\ldots,n\} \cong \A^1_{n*}$. 

(2)\ Suppose $\beta=0$. Then $x=\alpha$ and $z$ is free if $\alpha=x+z\beta$ has a solution. Further, $0=\alpha+yf(\alpha)$, 
which implies $\alpha=0$ if $f(\alpha)=0$. But this is impossible. So, $p^{-1}(k,0)=\emptyset$ if $k \in \{1,2,\ldots,n\}$.
If $f(\alpha) \ne 0$, $y$ is solved as $y=-\alpha/f(\alpha)$. So, $p^{-1}(\alpha,0)=\{(\alpha, -\alpha/f(\alpha),z) \mid 
z \in \C\}\cong \A^1$. Suppose $\beta \in \{1,2,\ldots,n\}$. Then $z=(\alpha-x)/\beta$ and the equation $x+yf(x)=\beta$ is 
written as $(x-\beta)(yg(x)+1)=0$, where $f(x)=g(x)(x-\beta)$. If $x=\beta$ then $y$ is free. If $yg(x)+1=0$ then 
$y=-1/g(x)$. This implies that $p^{-1}(\alpha,\beta)=\A^1\cup \A^1_{(n-1)*}$ with two components meeting in a single point
$(\beta, -1/g(\beta), (\alpha-\beta)/\beta)$ in the $(x,y,z)$-coordinates. By the above reasoning, we know that 
$\im(p)=\A^2\setminus S$ as stated above. The assertions (3) and (4) are also shown in the above argument.
\QED

If we restrict the morphism $p : \A^3 \to \A^2$ in the above proposition to a general linear plane in $\A^3$, we obtain 
an example of an endomorphism $\varphi : \A^2 \to \A^2$ which is not surjective. Such  examples have been constructed by 
Jelonek \cite{Jelonek}. 
\svskip

We further study the singular fibers of $\A^1_*$-fibrations.

\begin{example}\label{Example 2.17}{\em
Let $\wt{B}=\C[x,y,z,z^{-1}]$ be the Laurent polynomial ring in $z$ over a polynomial ring $\wt{A}=\C[x,y]$. Let $G_m$ 
and a cyclic group $G$ of order $n$ act on $\wt{B}$ by
\[
{}^t(x,y,z)=(x,y,tz)\quad \text{and}\quad {}^\zeta(x,y,z)=(x,\zeta y,\zeta^dz),
\]
where $\zeta$ is a primitive $n$-th root of $1$, $G$ is identified with the subgroup $\langle\zeta\rangle$ of $G_m$ 
generated by $\zeta$ and $d$ is a positive integer with $\gcd(d,n)=1$. Then the following assertions hold.
\begin{enumerate}
\item[(1)]
Let $B$ (resp. $A$) be the $G$-invariant subring of $\wt{B}$ (resp. $\wt{A}$) and let $d'$ be a positive integer such 
that $dd' \equiv 1 \pmod{n}$. Then $B=\C[x,y^n, y/z^{d'},z^n,z^{-n}]=\C[x,\eta,U,T,T^{-1}]/(U^n=\eta T^{-d'})$ and 
$A=\C[x,\eta]$, where $\eta=y^n,~U=y/z^{d'}$ and $T=z^n$.
\item[(2)]
Let $Y=\Spec B,~ X=\Spec(A)$ and $q : Y \to X$ be the morphism induced by the inclusion $A \hookrightarrow B$. Then 
$q$ is the quotient morphism of $Y$ by the $G_m$-action given by ${}^t(x,\eta,U,T)=(x,\eta,t^{-d'}U, t^nT)$. 
\item[(3)]
The singular locus $\Sing(q)$ is the line $\eta=0$ on $X$, and the singular fibers over $\Sing(q)$ are of the form 
$n\A^1_*$.
\end{enumerate} \QED }
\end{example}

In the above example, the $G$-action on $\wt{B}$ commutes with the $G_m$-action. Hence the $G_m$-action descends 
onto $Y$. Write $dd'=1+cn$. Then ${}^t(U^dT^c)=t^{-1}(U^dT^c)$. Hence the isotropy subgroup is trivial 
(resp. $G$) if $U \ne 0$ (resp. $U=0$). 

\begin{lem}\label{Lemma 2.18}
Let $p : Y \to X$ be an $\A^1_*$-fibration satisfying the following conditions.
\begin{enumerate}
\item[(1)]
$Y$ and $X$ are smooth, $p$ is equi-dimensional and $\Sing(p)$ is a smooth irreducible curve, say $C$. 
\item[(2)]
There exists a positive integer $n > 1$ such that the fiber $F_P$ over every point $P \in C$ is of the form $n\A^1_*$.
\item[(3)]
Either $X$ is factorial, or there exists an invertible sheaf $\SL$ on $X$ such that $\SL^{\otimes n} \cong \SO_X(-C)$. 
\end{enumerate}
Then there exists a cyclic covering $\mu : \wt{X} \to X$ of order $n$ ramifying totally over $C$ such that the 
normalization $\wt{Y}$ of $Y\times _X\wt{X}$ is an $\A^1_*$-bundle over $\wt{X}$, where $\nu : \wt{Y} \to Y\times_X\wt{X}$ 
is the normalization morphism. Further, there exists a $G_m$-action on $\wt{Y}$ such that the projection 
$\wt{q} : \wt{Y} \st{\mu} Y\times_X\wt{X} \to \wt{X}$ is the quotient morphism by the $G_m$-action. The quotient morphism 
$\wt{q}$ commutes with the cyclic covering group $G$, and hence descends down to the quotient morphism $q : Y \to X$ 
by the induced $G_m$-action which coincides with the given $\A^1_*$-fibration.
\end{lem}
\Proof
Write $X=\Spec A$. If $X$ is factorial, let $f \in A$ define the curve $C$ and let $\wt{X}=\Spec A[\xi]/(\xi^n-f)$.
If $\SO_X(-C)\cong \SL^{\otimes n}$, let $\wt{X}=\CSpec_X\bigoplus_{i=0}^{n-1}\SL^{\otimes i}$. Then the morphism 
$\mu : \wt{X} \to X$ is a cyclic covering of order $n$ totally ramifying over the curve $C$.

Let $\SU=\{U_i\}_{i \in I}$ be an affine open covering of $p^{-1}(C)$ and let $\eta_i$ be an element of 
$\Gamma(U_i,\SO_Y)$ such that $\eta_i=0$ defines $p^{-1}(C)_\red|_{U_i}$ for each $i$. We consider the case where 
$C$ is defined by $f=0$. The other case can be treated in a similar fashion. We can write $f=u_i\eta_i^n$ for 
$u_i \in \Gamma(U_i,\SO_Y^*)$. Then $\wt{Y}$ over $U_i$ is defined as $\CSpec_{U_i}\SO_{U_i}[\xi/\eta_i]$, where 
$(\xi/\eta_i)^n=u_i$. This implies that the morphism $\rho : \wt{Y} \st{\nu} Y\times_X\wt{X} \st{p_1} Y$ with the first
projection $p_1$ is a finite \'etale morphism, whence $\wt{Y}$ is smooth. It is clear that $\wt{q} : \wt{Y} \to 
\wt{X}$ is an $\A^1_*$-fibration. Since $\rho$ is finite and \'etale, $\wt{q}^{-1}(\wt{Q})$ with $\wt{Q} \in 
\mu^{-1}(C)_\red$ is a reduced $\A^1_*$ or a disjoint union of reduced $\A^1_*$'s. Let $D$ be a general smooth curve 
on $\wt{X}$ passing through $\wt{Q}$ and let $\wt{Z}$ be the normalization of $\wt{q}^{-1}(D)$. Then the canonical morphism 
$\ol{q} : \wt{Z} \to D$ is an $\A^1_*$-fibration on a normal affine surface $\wt{Z}$ and the fiber $\ol{q}^{-1}(\wt{Q})$ maps 
surjectively onto the fiber $\wt{q}^{-1}(\wt{Q})$. By Lemma \ref{Lemma 2.2}, $\ol{q}^{-1}(\wt{D})$ does not have two or 
more irreducible components which surject onto $\A^1_*$. Hence $\wt{q}^{-1}(\wt{Q})$ consists of a single reduced 
component isomorphic to $\A^1_*$. This implies that $\wt{Y}$ is an $\A^1_*$-bundle over $\wt{X}$.
  
An $\A^1_*$-bundle has a standard $G_m$-action along the fibers. Namely, if $\wt{q}^{-1}(V_j)\cong V_j\times \Spec \C[\tau_j,
\tau_j^{-1}]$ for an open covering $\SV=\{V_j\}_{j\in J}$ of $\wt{X}$, then $G_m$ acts as ${}^t\tau_j=t\tau_j$. This 
action clearly commutes with the action of the cyclic covering group $G$. Note that $\wt{Y}/G \cong Y$. Hence the 
$G_m$-action on $\wt{Y}$ descends down to $Y$, and gives rise to a $G_m$-action on $Y$. The quotient morphism 
$q : Y \to Y\quot G_m$ coincides with $p : Y \to X$ because $q=p$ on the open set $p^{-1}(X\setminus C)$.
\QED

A $G_m$-action of a smooth affine threefold is called {\em equi-dimensional} if the quotient morphism $q : Y \to X$ 
is equi-dimensional. 

\begin{lem}\label{Lemma 2.19}
Let $Y$ be a smooth affine threefold with a $G_m$-action and let $q : Y \to X$ be the quotient morphism with $X=Y\quot G_m$. 
Suppose that $q$ has equi-dimension one, $q$ is an $\A^1_*$-fibration and $X$ is smooth. Then a singular fiber of $q$ is 
either a tube or a cross. 
\end{lem}
\Proof
Let $F_P$ be a singular fiber $F$ of $q$. Let $C$ be a smooth irreducible curve on $X$ such that $P \in C$. Let 
$Z=Y\times_XC$ and let $\nu : \wt{Z} \to Z$ be the normalization morphism. Then the projection $\rho : Z \to C$ and 
$\wt{\rho}=\rho\cdot\nu : \wt{Z} \to C$ are the quotient morphisms by the induced $G_m$-actions on $Z$ and $\wt{Z}$. The 
fiber $\rho^{-1}(P)$ is the fiber $F_P$ and is the image of $\wt{F}_P=\wt{\rho}^{-1}(P)$ by $\nu$. By Lemma \ref{Lemma 2.2}, 
the fiber $\wt{F}_P$ is either a tube or a cross. If $\wt{F}_P$ is a tube, then $F_P$ is also a tube because 
$\nu$ restricted onto $\wt{F}_P$ commutes with the $G_m$-action. Suppose that $\wt{F}_P$ is a cross. Since the intersection 
point $\wt{Q}$ of two affine lines of $\wt{F}_P$ is a fixed point and $\nu_{\wt{F}_P}$ is surjective, either $F_P$ 
consists of two components meeting in a point $Q=\nu(\wt{Q})$ or $F_P$ is a contractible curve with $Q$ a fixed point and 
$F_P-\{Q\}$ a $G_m$-orbit. In the latter case, two branches of a cross are folded into a single curve. But this is impossible 
because the $G_m$-action on a cross viewed near the point $\wt{Q}$ has weights $-a, b$ respectively on two branches 
with $ab > 0$. So, $F_P$ consists of two branches meeting in one point $Q$. Consider the induced representation of $G_m$ on the 
tangent space $T_{Y,Q}$ at the fixed point $Q$. We can write it as ${}^t(x,y,z)=(x, t^{-a}y, t^bz)$ for a suitable system 
of local coordinates $\{x,y,z\}$. Then the fiber $F_P$ is given by $xy=0$ locally at $Q$. Hence two branches of $F_P$ 
meet transversally at $Q$, and $F_P$ is a cross.
\QED

We can prove a {\em converse} of this result.

\begin{thm}\label{Theorem 2.20}
Let $p : Y \to X$ be an $\A^1_*$-fibration on a smooth affine threefold $Y$. Suppose that $X$ is normal, $p$ is 
equi-dimensional and the singular fibers are tubes or crosses over the points of $\Sing(p)$ except for a finite 
set of points. Then there exists an equi-dimensional $G_m$-action on $Y$ such that the quotient morphism 
$q : Y \to Y\quot G_m$ coincides with the given $\A^1_*$-fibration $p$. Hence $X$ has at worst cyclic quotient singularities.
\end{thm}
\Proof
Since $\Sing(X)$ is a finite set, if a $G_m$-action is constructed on $p^{-1}(X\setminus\Sing(X))$ in such a way that 
the quotient morphism coincides with $p$ restricted on $p^{-1}(X\setminus\Sing(X))$, then the $G_m$-action extends to 
$Y$ by Hartog's theorem and the quotient morphism coincides with $p$. Since, as shown below, the construction of 
a $G_m$-action is local over $X\setminus\Sing(X)$, we may restrict ourselves to an affine open set of $X\setminus\Sing(X)$ 
and assume that $X$ is smooth from the beginning. 

By Lemmas \ref{Lemma 2.9} and \ref{Lemma 2.10}, the singular locus $\Sing(p)$ has no isolated points. Let 
$\Sing(p)=C_1\cup \cdots\cup C_r$ be the irreducible decomposition. Let 
\[
X_i=X\setminus(C_1\cup\cdots\cup \stackrel{\scriptsize \vee}{C_i}\cup\cdots \cup C_r)
\]
for $1 \le i \le r$, let $Y_i=p^{-1}(X_i)$ and let $p_i=p|_{Y_i} : Y_i \to X_i$. Then $Y_i$ and $X_i$ are affine and 
the $\A^1_*$-fibration $p_i$ has an irreducible singular locus, say $C_i$ by abuse of the notation. If $C_i$ has 
singular points, we replace $Y_i$ and $X_i$ first by $p_i^{-1}(X\setminus\Sing(C_i))$ and $X\setminus\Sing(C_i)$ 
respectively and then by $p_i^{-1}(U_{i\lambda})$ and $U_{i\lambda}$ respectively, 
where $\{U_{i\lambda}\}_{\lambda\in \Lambda}$ is an affine open covering of $X_i$. If there exist 
a family of equi-dimensional $G_m$-action $\sigma_i : G_m\times Y_i \to Y_i$ which induce the standard multiplication 
$\tau \mapsto t\tau$ for a variable on a general fiber $\tau$, then the actions $\{\sigma_i\}_{i=1}^r$ patch together 
and define an equi-dimensional $G_m$-action $\sigma : G_m\times Y \to X$. By Hartog's theorem, this $G_m$-action is 
also extended over the fibers lying over the singular points of $\Sing(p)$. In fact, the union of the fibers over the 
singular points of $\Sing(p)$ has codimension greater than one. 
   
So, we assume that $\Sing(p)$ is an irreducible smooth curve $C$. By shrinking $X$ again to a smaller open set, we may 
assume that the curve $C$ is principal, i.e., defined by a single equation $f=0$. If the singular fibers over 
the points of $C$ are tubes, the existence of a $G_m$-action follows from Lemma \ref{Lemma 2.18}. Suppose that the 
fibers over $C$ are crosses. Then $p^{-1}(C)=E_1\cup E_2$ with irreducible component $E_1, E_2$. Then $Y\setminus E_1$ 
and $Y\setminus E_2$ are affine and the fibration $p$ restricted to $Y\setminus E_1$ and $Y\setminus E_2$ have 
tubes over the curve $C$. Hence Lemma \ref{Lemma 2.18} implies that there exist $G_m$-actions on $Y\setminus E_1$ 
and $Y\setminus E_2$ and that they coincides on the general fibers of $p$. Hence they patch together and define a 
$G_m$-action on $Y$. 

The last assertion follows from Lemma \ref{Lemma 2.13}.
\QED

For examples with tubes or crosses as the singular fibers, we refer to Example \ref{Example 2.17} for tubes and Lemmas 
\ref{Lemma 4.1} and \ref{Lemma 4.4} below for crosses. In view of Lemma \ref{Lemma 2.2}, we have a satisfactory description 
on singular fibers of $\A^1_*$-fibrations on normal affine surfaces. If $Z$ is a normal affine surface with an 
$\A^1_*$-fibration $\rho : Z \to C$, then $p=\rho\times C' : Z \times C' \to C\times C'$ gives an $\A^1_*$-fibration 
on a normal affine threefold $Z\times C'$, where $C'$ is a smooth affine curve. Hence the same singular fibers as in the 
surface case appear in the product threefold case. But we can say much less in general. Let $p : Y \to X$ be an 
$\A^1_*$-fibration, where $Y$ is a smooth affine threefold and $X$ is a smooth affine surface. Let $F_P$ be the singular 
fiber of $p$ lying over a point $P \in X$. Let $C$ be a smooth curve on $X$ through $P$ and let $Z$ be the normalization 
of $Y\times_XC$. Then the induced morphism $p_C : Z \to C$ is an $\A^1_*$-fibration. Hence the fiber $F_P$ has 
a finite covering $\wt{F}_P \to F_P$, where $\wt{F}_P$ is a fiber of $p_C$ and hence has the form as described in 
Lemma \ref{Lemma 2.2}. We do not know exactly what the singular fiber $F_P$ itself looks like. 

Concerning the coexistence of tubes and crosses in the quotient morphism $q : Y \to X$ by a $G_m$-action, we have 
the following result.

\begin{lem}\label{Lemma 2.21}
Let $Y$ be a homology threefold with an effective, equi-dimensional $G_m$-action and let $q : Y \to X$ be the quotient 
morphism. Suppose that $q$ is an $\A^1_*$-fibration, $X$ is smooth and $\dim Y^{G_m}=1$. Then there are no tubes 
as singular fibers of $q$.
\end{lem}
\Proof
By the assertion (3) of Theorem \ref{Theorem 2.6}, $Y^{G_m} \cong \A^1$ and $q |_{Y^{G_m}} : Y^{G_m} \to X$ is a closed 
embedding. By Lemma \ref{Lemma 2.19}, any fiber through a point of $Y^{G_m}$ is a cross. Let $F=F_P$ be a tube of 
multiplicity $m > 1$. Let $p$ be a prime factor of $m$ and let $\Gamma = \Z/p\Z$. Consider the induced $\Gamma$-action 
on $Y$ as $\Gamma$ is a subgroup of $G_m$. Then $F$ is contained in the $\Gamma$-fixed point locus $Y^\Gamma$ (see the 
argument after the proof of Lemma \ref{Lemma 2.13}). Let $Q$ be a point on $F$. Then the induced tangential representation 
of $\Gamma$ on $T_{Y,Q}$ is written as 
\[
{}^\zeta(x,y,z)=(x,\zeta^a y,\zeta^b z), ~~ 0 \le a < p,~1 \le b <p,
\]
where $\zeta$ is a primitive $p$-th root of unity and $x$ is a coordinate with the tangential direction of the fiber. 
If $a> 0,~ b>0$, then the point $P=q(Q)$ is a singular point of $X$ (see the argument after the proof of 
Lemma \ref{Lemma 2.13}). Hence $a=0$. This implies that the component of $Y^\Gamma$ containing $F$ has dimension two. 
By the Smith theory, the locus $Y^\Gamma$ is a connected closed set. Hence there is a point $Q'$ of $Y^{G_m}$ such that 
the fiber $F'_{P'}$ through $Q'$ is contained in $Y^\Gamma$. Then we can write $F'_{P'}=m_1\A^1+m_2\A^2$ with $p \mid m_i$ 
for $i=1,2$.  But this is impossible. In fact, let $C$ be a smooth curve on $X$ passing through the point $P'$ and set 
$Z=Y\times_XC$. Then $Z$ has an induced $G_m$-action. Let $\wt{Z}$ be the normalization of $Z\times_C\wt{C}$, where 
$\wt{C} \to C$ is a finite covering of order $p$ totally ramifying over the point $P'$. Then there is an induced 
$G_m$-action on $\wt{Z}$ such that the induced morphism $\wt{q} : \wt{Z} \to \wt{C}$ is the quotient morphism. 
The fiber $\wt{q}^{-1}(\wt{P}')$ with $\wt{P}'$ a point of $\wt{C}$ lying over $P'$ consists of $p$-copies of the cross 
$\A^1+\A^1$ when taken with reduced structure. But this is impossible. 
\QED

Finally in this section, we discuss Question \ref{Question 2.12}. In fact, we prove a more general result.

\begin{lem}\label{Lemma 2.22}
Let $p : Y \to X$ be a dominant morphism from an algebraic variety of dimension $n+r$ to an algebraic variety 
$X$ of dimension $n$. Let $P$ be a point of $X$ and let $Q$ be a point of $Y$ such that $p(Q)=P$. Suppose that $Y$ is 
smooth at $Q$ and  the fiber $F:=p^{-1}(P)$ is reduced, $r$-dimensional and smooth at $Q$. Then $P$ is smooth at $X$.
\end{lem}
\Proof
Let $(R,\gm)$ (resp. $(S,\GM)$) be the local ring of $X$ (resp. $Y$) at $P$ (resp. $Q$). Since the fiber $F$ is 
defined by $\gm S$ at $Q$ and since $S/\gm S$ is a regular local ring as $Q$ is a smooth point of $F$, it follows 
that $\GM=\gm S+(z_{n+1},\ldots,z_{n+r})S$ for elements $z_{n+1}, \ldots, z_{n+r} \in \GM$. Since $(S,\GM)$ is 
a regular local ring of dimension $n+r$, we find $n$ elements $z'_1,\ldots,z'_n$ of $\gm$ with images $z_1, \ldots, z_n$ 
in $S$ such that $\{z_1,\ldots,z_n,z_{n+1},\ldots,z_{n+r}\}$ is a regular system of parameters of $(S,\GM)$. Since 
the completion $(\wh{S},\wh{\GM})$ is isomorphic to $\C[[z_1,\ldots,z_n,z_{n+1}, \ldots,z_{n+r}]]$, we can express 
any $h \in \wh{S}$ as a formal power series in $z_{n+1}, \ldots, z_{n+r}$ with coefficients in $\C[[z_1,\ldots,z_n]]$,
\[
h=\sum_{\vec{i}=0}^\infty \alpha_{\vec{i}}z^{\vec{i}},~~\alpha_{\vec{i}} \in \C[[z_1,\ldots,z_n]],
\]
where $z^{\vec{i}}= z_{h+1}^{i_{n+1}}\cdots z_{n+r}^{i_{n+r}}$ for $\vec{i}=(i_{n+1},\ldots,i_{n+r})$. We shall show that 
the completion $(\wh{R},\wh{\gm})$ is isomorphic to the formal power series ring $\C[[z_1,\ldots,z_n]]$. We consider 
the local complete intersection 
\[
H=\{z_j=0 \mid n+1 \le j \le n+r\}
\]
in $Y$ near the point $Q$ as a section transversal to the fiber $F$ at the point $Q$. We may assume that the restriction 
$p|_H : H \to X$ is quasi-finite. Hence the injective local ring homomorphism $R \to S \to S/zS$ induces an injective 
local homomorphism 
\[
(\wh{R},\wh{\gm}) \to (\wh{S},\wh{\GM}) \to (\wh{S}/z\wh{S},\wh{\GM}/z\wh{S}),
\]
because $(\wh{S}/(z_{n+1},\ldots,z_{n+r})\wh{S},\wh{\GM}/(z_{n+1},\ldots,z_{n+r})\wh{S})$ is viewed as the completion 
of the local ring $\SO_{H,Q}$. Sending $z_i$ to $z'_i,~i=1, \ldots, n$, we obtain a homomorphism 
\[
\C[[z_1,\ldots,z_n]] \lto \wh{R},
\]
which gives an isomorphism when composed with the mapping $\wh{R} \to (\wh{S}/z\wh{S})$. Let $h'$ be an element of $\wh{R}$ 
and $h$ its image in $\wh{S}$. With $h$ expressed as a power series as above, we find that $h'$ and the image in $\wh{R}$ 
of $\alpha_{\vec{0}} \in \C[[z_1,\ldots,z_n]]$ have the same image in $\wh{\GM}/z\wh{S}$ and hence coincide. Hence $\wh{R} 
\cong \C[[z_1,\ldots,z_n]]$. This implies that $X$ is smooth at $P$. 
\QED

\section{Homology threefolds with $\A^1$-fibrations}

In \cite{GMM}, it is shown that an $\A^1$-fibration $p : Y \to X$ from a smooth affine threefold to a normal surface has 
a factorization 
\[
p : Y \st{q} \wt{X} \st{\sigma} X,
\]
where $q : Y \to \wt{X}$ is the quotient morphism by a $G_a$-action and $\sigma : \wt{X} \to X$ is the birational morphism 
such that $\wt{A}=\Gamma(\wt{X},\SO_{\wt{X}})$ is the factorial closure of $A=\Gamma(X,\SO_X)$ in $B=\Gamma(Y,\SO_Y)$. 
Namely, 
\[
\wt{A}=\{b \in B \mid \text{$b$ is a factor of an element $a \in A$}\}. 
\]
Thus we may look into the quotient morphism $q : Y \to X$ by a $G_a$-action. We set 
$\Sing(q)=\{P \in X \mid F_P \not\cong \A^1\}$, where $F_P$ is the fiber over $P \in X$, and call it the {\em singular locus} 
of $q$. 

\begin{lem}\label{Lemma 3.1}
Let $q : Y \to X$ be the quotient morphism of a smooth affine threefold $Y$ with respect to a $G_a$-action. Suppose that 
$q$ is equi-dimensional. Then the following assertions hold.
\begin{enumerate}
\item[(1)]
If a fiber $F_P$ has an irreducible component which is reduced in $F_P$, then the point $P$ is smooth in $X$.
\item[(2)]
$\Sing(q)$ is a closed set.
\end{enumerate}
\end{lem}
\Proof
(1)\ The assertion follows from Lemma \ref{Lemma 2.22}.

(2)\ Since $q$ has equi-dimension one, the fixed point locus $Y^{G_a}$ consists of fiber components \cite[Corollary 3.2]{GMM}. 
Then every fiber of $q$ is a disjoint union of contractible curves \cite[Lemma 3.5]{GMM}. Further, a contractible irreducible 
component is isomorphic to $\A^1$ if it contains a non-fixed point, or is contained in $Y^{G_a}$. 

It suffices to show that given a fiber $F_P$ isomorphic to $\A^1$ there exists an open neighborhood $U$ of $P$ in $X$ 
with $F_{P'} \cong \A^1$ for all $P' \in U$. The subsequent proof is similar to the one for $G_m$-actions 
(cf. Lemma \ref{Lemma 2.11}). Let $C$ be an irreducible curve on $X$ through $P$ and let $Z$ be 
the normalization of $Y\times_XC$. Suppose that the fiber of $q$ over a general point of $C$ consists of $m$ copies of 
$\A^1$, where $m > 1$. Then $Z$ has an induced $G_a$-action such that 
\[
\CD
Z~~ @>\mu >> Y\times_XC \\
@VV\wt{q} V                  @VV q V \\
\wt{C}~~  @>\nu >> C
\endCD
\]
where $\wt{q}$ is the quotient morphism and $\wt{C}$ is the normalization of $C$ in $Z$. The morphism 
$\nu : \wt{C} \to C$ is a finite covering of degree $m$ (the Stein factorization), which is ramified over the point $P$. 
It then follows that the fiber $F_P$ is non-reduced, and this is a contradiction.  

Consider the closure $\ol{\Sing(q)}$ and remove from $X$ all the irreducible components of $\ol{\Sing(q)}$ over which 
a general fiber of $q$ consists of more than one irreducible component. Thus we may assume that all singular fibers of 
$q$ are of the form $m\A^1$ with $m > 1$. Note that $P$ is a smooth point by the assertion (1). Replacing $X$ by an 
affine open neighborhood $U$ of $P$ and accordingly $Y$ by the inverse image $q^{-1}(U)$, we may assume that $X$ is smooth. 
Consider the tangential homomorphism of the tangent bundles 
\[
dq : \ST_Y \to q^*\ST_X
\]
and let $\SC$ be the cokernel of $dq$. Then $\SC$ is a coherent $\SO_Y$-module. The support $T=\Supp(\SC)$ is a closed set 
such that $T=q^{-1}(q(T))$. If $P \in \ol{\Sing(q)}$, then $F_P \subset T$, but $F_P \cap T=\emptyset$. This is a 
contradiction. Then, by Dutta \cite[Theorem]{Dutta}, there exists an open neighborhood $U$ of $P$ such that $q^{-1}(U)$ 
is an $\A^1$-bundle over $U$. Hence $\Sing(q)$ is a closed set.
\QED

Concerning a $G_a$-action, we ask the following

\begin{question}\label{Question 3.2}{\em
Let $Y$ be an affine variety with a $G_a$-action. Suppose that the algebraic quotient $X=Y\quot G_a$ exists as an affine 
variety, i.e., the $G_a$-invariant subring of $\Gamma(Y,\SO_Y)$ is an affine domain. Let $Q$ be a point of $Y$ with trivial 
isotropy group, i.e., a point which is not $G_a$-fixed. Is then the fiber $F_P$ of the quotient morphism $q : Y \to X$ 
passing through the point $Q$ reduced, where $P=q(Q)$ ?}
\end{question}

The answer is negative and an example is given by an affine pseudo-plane \cite[Corollary 3.16]{GMM}. See also \cite[Lemma 2.4]{Mi7}.

\begin{example}\label{Example 3.3}{\em Let $C$ be the smooth conic in $\BP^2=\Proj~\C[X_0,X_1,X_2]$ defined by $X_2^2=X_0X_1$. 
Let $L$ be the tangent line of $C$ at the point $(X_0,X_1,X_2)=(0,1,0)$. Hence $L$ is defined by $X_0=0$. Let $\Lambda$ 
be the linear pencil generated by $C$ and $2L$. Then the $G_a$-action on $\BP^2$ defined by 
\[
{}^tX_0=X_0,~~{}^tX_1=X_1+2X_2 t+X_0t^2,~~{}^tX_2=X_2+X_0t
\]
has a unique fixed point $Q$ and preserves the members of the pencil $\Lambda$. Let $X$ be the complement of $C$ in $\BP^2$.
Then $\Lambda$ defines an $\A^1$-fibration $q : X \to \A^1$ which turns out to be the quotient morphism of the induced 
$G_a$-action on $X$. Although the $G_a$-action has no fixed points, $q$ has a multiple fiber $2\ell$ which comes from the 
member $2L$ of $\Lambda$. The affine surface $X$ is an affine pseudo-plane and its universal covering is a Danielewski 
surface.}
\end{example}

This is the case for every affine pseudo-plane $q : X \to Z$ of type $(m,r)$ with an integer $r \ge 1$ such that  
$r \equiv 1 \pmod{m}$. For the definition of an affine pseudo-plane of type $(m,r)$, see \cite{MM} where the type is 
denoted by $(d,r)$ instead of $(m,r)$. In particular, $q : X \to Z$ is an $\A^1$-fibration with a unique singular fiber 
$F_0$ of multiplicity $m > 1$, i.e., $F_0=m\A^1$ and $Z \cong \A^1$. $X$ is denoted by $X(m,r)$.

\begin{lem}\label{Lemma 3.4}
Let $q : X \to Z$ be an affine pseudo-plane of type $(m,r)$ with $r \equiv 1 \pmod{m}$. Then $q$ is given by a fixed-point 
free $G_a$-action.
\end{lem}
\Proof
Let $\wt{X}(m,r)$ be the universal covering of $X(m,r)$. Let $H(m)=\Z/m\Z$ be the covering group which is identified with 
the $m$-th roots of unity. By \cite[Lemma 2.6]{MM}, $\wt{X}(m,r)$ is isomorphic to a hypersurface in $\A^3=\Spec \C[x,y,z]$ defined by 
\[
x^rz+(y^m+a_1xy^{m-1}+\cdots+ a_{m-1}x^{m-1}y+a_mx^m)=1,~~ a_i \in \C.
\]
The Galois group $H(m)$ acts as 
\[
{}^\lambda(x,y,z)=(\lambda x,\lambda y, \lambda^{-r}z), ~~\lambda \in H(m).
\]
The projection $(x,y,z) \mapsto x$ defines an $\A^1$-fibration $\wt{q} : \wt{X}(m,r) \to \wt{Z}$, where $\wt{Z}\cong \A^1$. 
Further, there is a $G_a$-action on $\wt{X}(m,r)$ defined by 
\begin{eqnarray*}
\lefteqn{{}^t(x,y,z)=} \\
&& (x,y+tx^r, z-x^{-r}\{((y+tx^r)^m+a_1x(y+tx^r)^{m-1}+\cdots+a_mx^m) \\
&& \qquad \qquad -(y^m+a_1xy^{m-1}+\cdots+a_{m-1}x^{m-1}y+a_mx^m)\}).
\end{eqnarray*}
Then it follows that 
\begin{enumerate}
\item[(1)]
the fiber $\wt{F}_0$ of $\wt{q}$ over the point $x=0$ is a disjoint union of $m$-copies of the affine line,
\item[(2)]
the $G_a$-action preserves the fibration $\wt{q}$,
\item[(3)]
the $G_a$-action preserves and acts non-trivially on each connected component of $\wt{F}_0$, 
\item[(4)]
if $r \equiv 1 \pmod{m}$ then the $H(m)$-action commutes with the $G_a$-action.
\end{enumerate}
Hence the $G_a$-action descends to a $G_a$-action on $X(m,r)$ which has no fixed points.
\QED

We consider, however, a result implied by the assumptions in Question \ref{Question 3.2}. We need some auxiliarly results. 
Let $B=\Gamma(Y,\SO_Y)$ and let $A=\Gamma(X,\SO_X)$. Let $\delta$ be a locally nilpotent derivation on $B$ which 
corresponds to the given $G_a$-action. Further, let $\gm$ be the maximal ideal of $A$ corresponding to the point $P$. 
First we recall the following result in \cite[Theorem 3.3]{Mi7}.
  
\begin{lem}\label{Lemma 3.5}
Let $B$ be an affine domain over $\C$ and let $\delta$ be a locally nilpotent derivation on $B$. Let $A = \Ker \delta$.
Suppose that $B/\gm B$ is an integral domain over $\C$ of dimension one and that the associated $G_a$-action on $\Spec B/\gm B$ 
has no fixed points. Then the following assertions hold.
\begin{enumerate}
\item[(1)] 
For any integer $n > 0$, there exists an element $z_n \in B$ such that $B/\gm^n B = R_n[z_n]$, where $R_n$ is an Artin local ring 
and we denote the residue class of $z_n$ in $B/\gm^n B$ by the same letter. 
\item[(2)] 
For $m > n$, we have a natural exact sequence
\[
0 \to \gm^n R_m \to R_m \st{\theta_{nm}} R_n \to 0,
\]
where $\theta_{nm}$ is a local homomorphism. Let $\wh{R}=\varprojlim_n R_n$. Then $\wh{R}$ is a complete local ring.
\end{enumerate}
\end{lem}

Let $\SO$ be the local ring $\SO_{X,P}$ and denote the ideal $\gm\SO$ by the same letter $\gm$. Let $\wh{\SO}$ 
be the $\gm$-adic completion of $\SO$. Then $B\otimes_A\SO$ is an affine domain over $\SO$ with the associated 
locally nilpotent derivation $\delta$ such that $\delta$ is nonzero and $\Ker\delta=\SO$. Further {\em we assume the following 
condition {\rm (H)} that $(\gm^n(B\otimes_A\SO))\cap A=\gm^n$ for every $n > 0$.} This condition is satisfied if $B\otimes_A\SO$ 
is $\SO$-flat (see \cite[Lemma 1.4]{GM2}).

\begin{prop}\label{Proposition 3.6}
With the above notations and assumptions, the following assertions hold.
\begin{enumerate}
\item[(1)]
$(B\otimes_A\SO)/\gm^n(B\otimes_A\SO)\cong B/\gm^n B=R_n[z_n]$.
\item[(2)]
The projective limit $\wh{B}=\varprojlim_n (B\otimes_A\SO)/\gm^n(B\otimes_A\SO)$ contains $B\otimes_A\SO$ as a subring. 
Further, $\wh{B}$ has a derivation $\wh{\delta}$ and an element $\wh{z}$ such that $\wh{\delta}\mid_{B\otimes_A\SO}=\delta$,
$\Ker\wh{\delta}=\wh{R}$ and $\wh{\delta}(\wh{z})=1$. The ring $\wh{B}$ itself is a subring of a formal power series ring 
$\wh{R}[[\wh{z}]]$.
\item[(3)]
Let $\wh{B}_{\rm ln}$ be the subring of $\wh{B}$ consisting of elements for which $\wh{\delta}$ is locally nilpotent. 
Then $\wh{B}_{\rm ln}=\wh{R}[\wh{z}]$ and $B\otimes_A\SO$ is a subring of $\wh{B}_{\rm ln}$.
\end{enumerate}
\end{prop}
\Proof
(1)\ The assertion follows from Lemma \ref{Lemma 3.5}.

(2)\ It suffices to show that if $b$ is an element of $\bigcap_n \gm^n(B\otimes_A\SO)$ then $b=0$. Suppose that 
$b \ne 0$. Then there exists an integer $r \ge 0$ such that $\delta^r(b) \ne 0$ and $\delta^{r+1}(b)=0$. Then 
$\delta^r(b) \in \Ker \delta=\SO$. Since $\delta(\gm^n(B\otimes_A\SO)) \subseteq \gm^n(B\otimes_A\SO)$ for every $n > 0$, 
by the condition (H), we have $\delta^r(b) \in (\bigcap_n \gm^n(B\otimes_A\SO))\cap \SO =\bigcap_n \gm^n$, which is 
zero by Krull's intersection theorem. Hence $\delta^r(b)=0$. This is a contradiction. Hence $B\otimes_A\SO$ is a subring 
of $\wh{B}$. The sequence $\{z_n\}_{n \ge 1}$ is a Cauchy sequence because $z_m-z_n \in \gm^n(B\otimes_A\SO)$ for $m > n$. 
Hence it converges to an element $\wh{z}$ of $\wh{B}$. The derivation $\delta$ extends to a derivation $\wh{\delta}$. 
Since we may assume that $\delta(z_n)\equiv 1 \pmod{\gm^n(B\otimes_A\SO)}$, it follows that $\wh{\delta}(\wh{z})=1$ 
in $\wh{B}$. Let $b_n$ be an element of $B\otimes_A\SO$. Since $(B\otimes_A\SO)/\gm^n(B\otimes_A\SO) \cong R_n[z_n]$ 
and $\wh{z}-z_n \in \gm^n\wh{B}$, there exists a polynomial $f_n(\wh{z}) \in \wh{R}[\wh{z}]$ such that 
$b_n-f_n(\wh{z}) \in \gm^n\wh{B}$. This implies that a Cauchy sequence $\{b_n\}_{n \ge 0}$ in $B\otimes_A\SO$ is 
approximated by a Cauchy sequence $\{f_n(\wh{z})\}_{n\ge 0}$ in $\wh{B}$. So, $\wh{B}$ is a subring of $\wh{R}[[\wh{z}]]$.

(3)\ Let $b$ be an element of $B\otimes_A\SO$. Write 
\[
b=\alpha_0+\alpha_1\wh{z}+\cdots+\alpha_i\wh{z}^i+ \cdots~=~\sum_{i=0}^\infty \alpha_i\wh{z}^i
\]
as an element of $\wh{R}[[\wh{z}]]$. If this is an infinite series, the module $M(b)$ generated over $\wh{\SO}$ by 
$\wh{\delta}^i(b)$ for $i \ge 0$ is not finitely generated over $\wh{\SO}$ because it is not finitely generated over $\wh{R}$.
Meanwhile, $\wh{\delta}^i(b)=\delta^i(b)$ for every $i \ge 0$ and $M(b)$ is finitely generated over $\wh{\SO}$ because 
$\delta$ is locally nilpotent. Hence $b \in \wh{R}[\wh{z}]$. This implies that $B\otimes_A\SO$ is a subring of 
$\wh{R}[\wh{z}]$. It is clear that $\wh{B}_{\rm ln}=\wh{R}[\wh{z}]$.
\QED

We recall the following result of Kaliman-Saveliev \cite[Corollary 2.8]{KS}. 

\begin{lem}\label{Lemma 3.7}
Let $Y$ be a smooth contractible affine threefold with a fixed-point free $G_a$-action. Then $Y$ is isomorphic to 
$X\times \A^1$ with $G_a$ acting on the second factor, where $X=Y\quot G_a$.
\end{lem}

This result inspires us a challenging problem.

\begin{question}\label{Question 3.8}{\em Let $Y$ be a homology threefold with a fixed-point free $G_a$-action. Is $Y$ 
$G_a$-equivariantly isomorphic to a product $X \times \A^1$ with $G_a$ acting on the second factor, where $X=Y\quot G_a$ ?}
\end{question}

If $Y$ is replaced by a $\Q$-homology threefold, the answer is negative.

\begin{remark}\label{Remark 3.9}{\em Let $q : X \to Z$ be an affine pseudo-plane as in Lemma \ref{Lemma 3.4}. Then 
$X$ is a $\Q$-homology plane. Let $Y=X \times \A^1$ with $G_a$ acting on $X$ and trivially on the second factor. Then 
$Y$ is a $\Q$-homology threefold, $Y \quot G_a=Z \times \A^1$ and the quotient morphism is 
$q_Y=q \times \id_{\A^1} : X\times \A^1 \to Z\times \A^1$. Hence $\Sing(q_Y) \cong \A^1$.}
\end{remark}

\section{Contractible affine threefolds with $\A^1_*$-fibrations}

We first discuss Question \ref{Question 2.7}. 

\begin{lem}\label{Lemma 4.1}
Let $X$ be a homology plane of log Kodaira dimension $-\infty$ or $1$ and let $C$ be a curve on $X$ isomorphic 
to $\A^1$. Let $V=X\times \A^1$ with a $G_m$-action induced from the standard action on the $\A^1$-factor, i.e., ${}^t(Q,x)
=(Q,tx)$, where $Q \in X$ and $x$ is a coordinate of $\A^1$. Let $\sigma : W \to V$ be the blowing-up of $V$ with center 
$C\times (0)$ which is identified with $C$. Let $Y=W\setminus\sigma'(X\times (0)$, where $\sigma'(X\times (0))$ is the proper 
transform of $X\times (0)$. Then the following assertions hold.
\begin{enumerate}
\item[(1)]
$Y$ is a homology threefold with an induced $G_m$-action such that the quotient morphism $q : Y \to Y\quot G_m$ has 
the quotient space $Y\quot G_m \cong X$. If $X$ is contractible, the threefold $Y$ is also contractible.
\item[(2)]
The singular locus $\Sing(q)$ is the curve $C$. For every point $P$ of $C$, the fiber $q^{-1}(P)$ consists of two affine 
lines meeting in one point. Thus we have the situation treated in Theorem \ref{Theorem 2.6}, (3).
\item[(3)]
If $\lkd(X)=1$, then $\lkd(Y)=1$. Hence $Y$ cannot be written as $Y \cong Z \times Y^{G_m}$, i.e., the answer to the 
question \ref{Question 2.7} is negative.
\item[(4)]
If $\lkd(X)=-\infty$, i.e., $X \cong \A^2$, then $Y \cong Z\times Y^{G_m}$ with $Z$ a smooth affine surface with a 
$G_m$-action.
\end{enumerate}
\end{lem}
\Proof
(1)\ Set $D=X\times (0),~ L=C\times \A^1,~ L'=\sigma'(L)$ and $E'=\sigma^{-1}(C)\setminus\sigma'(D)$. Then $\sigma : 
(Y, E') \to (V,D)$ is an affine modification with $\sigma(E')=C$ (see \cite{KZ}). Write $X=\Spec A$. Then $A[x]$ 
is the coordinate ring of $V$. Since $A$ is factorial by Lemma \ref{Lemma 1.1}, the curve $C$ is defined by an element 
$\xi$ of $A$. The hypersurface $D=X\times (0)$ in $V$ is defined by $x=0$. Let $I$ be the ideal of $A[x]$ generated 
by $\xi$ and $x$. Then $Y$ has the coordinate ring $\Sigma_{I,x}(A[x])$ which is the affine modification of $A[x]$ 
along $(x)$ with center $I$. Clearly, $Y$ is a smooth affine threefold. By \cite[Proposition 3.1 and Theorem 3.1]{KZ}, 
it follows that $Y$ is a homology threefold and is contractible provided so is $X$. Let $q_0 : V \to X$ be the projection 
to $X$, which is in fact the quotient morphism by the $G_m$-action such that $V^{G_m}=D$. By the above process, the 
$G_m$-action is inherited on $Y$ and the quotient morphism $q : Y \to Y\quot G_m$ is induced by $q_0$. 

(2)\ Meanwhile, in passing from $V$ to $Y$, the fiber over a point $Q \in X\setminus C$ loses the point $(Q,0)$ in $D$ 
and becomes isomorphic to $\A^1_*$. The fiber $q^{-1}(Q)$ is a cross $\A^1+\A^1$ with two $\A^1$ meeting in the point 
$q^{-1}\cap(L'\cap E')$. Hence $E'\cap L' =Y^{G_m}$ and $\Sing(q)=C$. 

(3)\ Suppose that $\lkd(X)=1$. In fact, there is a unique affine line lying on $X$. Since the general fibers of $q$ are 
isomorphic to $\A^1_*$, we have an inequality $\lkd(Y) \ge \lkd(X)+\lkd(F)$ by Kawamata \cite{Ka2}, where $F$ is 
a general fiber of $q$. Hence $\lkd(Y) \ge 1$. Furthermore, $X$ itself has an $\A^1_*$-fibration $\pi : X \to B$ such that 
$C$ is a fiber of $\pi$ and $\pi^{-1}(U) \cong U\times \A^1_*$, where $U$ is an open set of $B$ contained in $\A^1_{2*}$.
By the above construction, $q^{-1}(\pi^{-1}(U)) \cong \pi^{-1}(U)\times \A^1_* \cong U\times \A^1_*\times\A^1_*$. 
Hence $\lkd(Y) \le \lkd(U\times \A^1_*\times\A^1_*) =1$. This implies that $\lkd(Y)=1$. If $Y \cong Z\times \A^1$ as 
inquired in Question \ref{Question 2.7}, then it would follow that $\lkd(Y)=-\infty$. But this is not the case.

(4)\ Suppose that $\lkd(X)=-\infty$. Then $X \cong \A^2$, and the affine line $C$ is chosen to be a coordinate line by 
AMS theorem. Namely, there exists a system of coordinates $(\xi,\eta)$ of $X$ such that $C$ is defined by $\xi=0$. The 
affine modification $\Sigma_{I,x}(A[x])$ is equal to $\C[x, \xi/x, \eta]$. Set $y=\xi/x$ and $R=\C[x,y]$. Then the 
induced $G_m$-action on $\Sigma_{I,x}(A[x])$ is given by ${}^t(x,y,\eta)=(tx,t^{-1}y,\eta)$. Hence the threefold $Y$ 
is isomorphic to $Z\times\A^1$, where $Z=\Spec R\cong \A^2$ and $\A^1=\Spec \C[\eta]$. So, the answer to Question 
\ref{Question 2.7} is affirmative.
\QED

\begin{remark}\label{Remark 4.2}{\em
The quotient morphism $q : Y \to X$ has crosses $\A^1+\A^1$ with each multiplicity one as the singular fibers over 
$\Sing(q)\cong \A^1$. The locus of intersection points of crosses is the fixed point locus $\Gamma=Y^{G_m}$. 
There are two embedded affine planes $Z_1, Z_2$ meeting transversally along $\Gamma$. If $Y$ is a homology threefold, 
$Z_1, Z_2$ are defined by $f_1=0, f_2=0$. Hence $\Gamma$ is defined by the ideal $I_1=(f_1,f_2)$ of $B_1=\Gamma(Y,\SO_Y)$. 
The affine modification $B_2=\Sigma_{I_1,f_1}(B)$ or $B_2'=\Sigma_{I,f_2}(B)$ gives rise to a smooth affine threefold 
$Y_2=\Spec B_2$ or $Y'_2=\Spec B_2'$ with an equi-dimensional $G_m$-action, which gives the quotient morphism 
$q_2 : Y_2 \to X$ or $q'_2 : Y'_2 \to X$. Both $Y_2$ and $Y'_2$ are homology threefolds (resp. contractible threefolds) 
if $Y$ is a homology threefold (resp. contractible threefold). A difference between $Y_2$ (or $Y'_2$) and $Y$ is 
that the crosses have multiplicities. Namely, $Y_2$ (resp. $Y'_2$) has crosses $2\A^1+\A^1$ (resp. $\A^1+2\A^1$). 
We can repeat this process to produce homology threefolds or contractible threefolds which have crosses with higher 
multiplicities. If a cross is written as $m\A^1+n\A^1$, then $\gcd(m,n)=1$ (see the argument in the last part of 
the proof of Lemma \ref{Lemma 2.21}).  }
\end{remark}

Another remark to Lemma \ref{Lemma 4.1} is the following. 

\begin{remark}\label{Remark 4.3}{\em 
Take a $\Q$-homology plane $X$ instead of a homology plane in the construction of $Y$ in Lemma \ref{Lemma 4.1}. 
Then we can consider such $\Q$-homology planes with log Kodaira dimension $-\infty, 0$ and $1$. We take an 
embedded line $C$ in $X$ and blow up the center $C\times (0)$ in $V=X\times\A^1$. By the same construction, we obtain 
a smooth affine threefold $Y$ with an equi-dimensional $G_m$-action. The threefold $Y$ has the same homology groups 
as $X$. Hence $Y$ is a $\Q$-homology threefold, which is not of the product type $Z\times \A^1$ with a 
$\Q$-homology plane $Z$ provided $\lkd(X) \ne -\infty$. For the embedded lines in the case of log Kodaira dimension $0$, 
see \cite[Theorem]{GP} for a complete classification.}
\end{remark}

In Lemma \ref{Lemma 4.1}, the center of blowing-up is $C\times (0)$, where $C$ is an embedded line in a homology plane,
and the resulting homology threefold has log Kodaira dimension at most one. We can apply a process similar to the one 
used in Lemma \ref{Lemma 4.1} with a point as the center to obtain the first assertion of the following lemma. This 
is first constructed in \cite[Example 3.7]{KR2}.

\begin{lem}\label{Lemma 4.4}
{\em (1)}\ \ Let $X$ be a homology plane and let $P_0$ be a point of $X$. Let $Y_0=X \times \A^1$ with $G_m$ 
acting trivially on $X$ and on $\A^1$ with weight $-1$. Let $Q_0=P_0\times (0)$. Blow up the point $Q_0$ and remove 
the proper transform of $X\times (0)$ to obtain a smooth affine threefold $Y$. Then $Y$ is a homology threefold 
with a $G_m$-action such that the quotient morphism $q : Y \to X$ is induced by the first projection 
$p_1 : Y_0 \to X$, the fixed point locus $Y^{G_m}$ consists of the unique point which is the intersection point 
of the proper transform of $P_0\times \A^1$ with the exceptional surface of the blowing-up (whence the $G_m$-action 
on $Y$ is {\em hyperbolic}) and $\lkd(Y)=\lkd(X)$. If $X$ is contractible, so is $Y$.

{\em (2)}\ \ Let $Y$ be a smooth homology threefold with a hyperbolic $G_m$-action. Let $Q_0$ be the unique fixed 
point and let $Z_0$ be the two-dimensional fiber component of the quotient morphism $q : Y \to X$, where $Q_0 \in Z_0$ 
and $X=Y\quot G_m$ (cf. Theorem \ref{Theorem 2.6}). Let $Y'$ be the affine modification $\Sigma_{Q_0,Z_0}(Y)$ 
which is the blowing-up of $Y$ at the center $Q_0$ with the proper transform of $Z_0$ deleted off. Then $Y'$ is 
a smooth homology threefold. If $Y$ is contactible, then so is $Y'$.
\end{lem}
The statement and the proof depend on \cite[\S 3]{KZ}.

\begin{example}\label{Example 4.5}{\em 
Let $Y$ be the Koras-Russell threefold $x+x^2y+z^2+t^3=0$. Then a hyperbolic $G_m$-action on $Y$ is given by 
\[
{}^\lambda(x,y,z,t)=(\lambda^6x,\lambda^{-6}y,\lambda^3z,\lambda^2t),~~\lambda \in \C^*.
\]
The fixed point $Q_0$ is $(0,0,0,0)$, and the two-dimensional fiber component $Z_0$ is defined by $y=0$. Let $B$ 
be the coordinate ring of $Y$ and let $I=(x,y,z,t)$ which is the maximal ideal of $Q_0$. The affine modification 
$B'=\Sigma_{I,y}(B)$, which is the coordinate ring of $Y'$, is given as
\[
B'=\C[x',y,z',t'],~~x'=\frac{x}{y},~z'=\frac{z}{y},~t'=\frac{t}{y}.
\]
Hence $Y'$ is a hypersurface $x'+{x'}^2y^2+{z'}^2y+{t'}^3y^2=0$ and the hyperbolic $G_m$-action is given by 
\[
{}^\lambda(x',y,z',t')=(\lambda^{12}x',\lambda^{-6}y,\lambda^9z',\lambda^8t'),~~\lambda \in \C^*.
\]
We can further repeat the affine modifications of the same kind to $Y', Y''$ etc.

The above Koras-Russell threefold and its affine modifications are examples of smooth contractible threefolds $Y$ 
with a hyperbolic $G_m$-action such that the quotient $X$ is isomorphic to that of the tangent space $T_{Q_0}$ at 
the unique fixed point $Q_0$ of $Y$ by the induced tangential representation. In \cite[Theorem 4.1]{KR2} a description of all 
such threefolds is given. 

In \cite{KML} the Makar-Limanov invariants $\ML(Y)$ of such threefolds are computed. 
To apply this result the equation for $Y'$ has to be brought into a standard form that exhibits $Y'$ as a cyclic 
cover of $\A^3$. To this end let $B''={B'}^{\omega}=\C[x',y,\zeta, t']$ with a square root $\omega$ of the unity and 
$\zeta ={z'}^2$. Put
\[
x''=-({x'}^2y+\zeta+{t'}^3y).
\]
Then $x'=x''y$ and we see that $B''=\C[x'', y, t']$ is a polynomial ring in three variables and 
\[
x''+{x''}^2y^3+\zeta+{t'}^3y=0.
\]
Hence it follows that $B'=\C[x'', y, z', t']$ with a defining equation 
\[
x''+{x''}^2y^3+{z'}^2+{t'}^3y=0.
\] 
It now follows from \cite[Theorem 8.4]{KML} that $\ML(Y)=\C[x]$ and $\ML(Y')$ $=B'$. \QED}
\end{example}

In order to show that a fixed point exists under a $G_m$-action on a $\Q$-homology threefold $Y$, we used 
the Smith theory and its variant. The following result without using the Smith theory is of some interest. 

\begin{lem}\label{Lemma 4.6}
Let $Y$ be a smooth affine variety with an effective $G_m$-action. Suppose that there are no fixed points. 
Then the Euler number $e(Y)$ of $Y$ is zero.
\end{lem}
\Proof
Let $q : Y \to X$ be the quotient morphism. Since there are no fixed points, every fiber is isomorphic to $\A^1_*$ 
when taken with reduced structure. The general fibers of $q$ are reduced $\A^1_*$ and special fibers are multiple 
$A^1_*$. We work with the complex analytic topology. Considering the isotropy groups of the fibers, there exists 
a descending chain of closed subsets of $X$
\[
F_0 \supset F_1 \supset \cdots \supset F_i \supset F_{i+1} \supset \cdots 
\]
such that $F_0=X$ and the isotropy group of the fiber over a point of $F_i-F_{i+1}$ is constant, say $G_i$. Then 
$F_i-F_{i+1}$ is covered by open sets $\{U_{i\lambda}\}_{\lambda\in \Lambda_i}$ such that $q^{-1}(U_{i\lambda})_\red
\cong G_m\times_{G_i}V_{i\lambda}$, where $V_{i\lambda}$ is a suitable slice (cf. \cite[Th\'eor\`eme du slice \'etale 
et Remarque $3^\circ$]{Luna}). Hence $q^{-1}(F_i)_\red$ 
is a $\C^*$-bundle over the open set $F_i-F_{i+1}$. This implies that the Euler number $e(Y)$ is zero. 
\QED

Lemma \ref{Lemma 4.6} implies that any non-trivial $G_m$-action on a smooth affine variety $Y$ has a fixed point 
if $e(Y)\ne 0$. Furthermore, considering the induced tangential representation at a fixed point, we know that 
the fixed point locus $Y^{G_m}$ is smooth. However, we do not know if the fixed point locus $Y^{G_m}$ is connected. 
The following example shows that the connectedness fails in general.

\begin{example}\label{Example 4.7}{\em
Let $X=\A^2$ and let $Y_0=X\times \A^1$ which has a standard $G_m$-action on the factor $\A^1$ with the point $(0)$ 
as a fixed point. Choose two parallel lines $\ell_1, \ell_2$ on $X$, and let $Z_i=\ell_i\times \A^1$ for $i=1,2$. 
Let $\sigma : W \to Y_0$ be the blowing-ups with centers $\ell_1\times (0)$ and $\ell_2\times (0)$ and let $Y$ be 
$W$ with the proper transform of $X\times (0)$ removed. Then the restriction of $p_1\cdot\sigma$ onto $Y$ gives 
a morphism $q : Y \to X$ which is eventually the quotient morphism by the $G_m$-action on $Y$ induced by the one 
on $Y_0$, where $p_1$ is the projection $Y_0 \to X$. The singular fibers of $q$ are crosses over $\ell_1\cup \ell_2$. 
The locus of intersection points of crosses is the fixed point locus $Y^{G_m}$. Hence $Y^{G_m}$ is not connected. 
The Euler number $e(Y)$ of the threefold $Y$ is two. \QED}
\end{example}

In Section 2, we observed a $G_m$-action on a $\Q$-homology threefold $Y$ such that the quotient morphism 
$q : Y \to X$ has relative dimension one. The following result deals with the case of $q$ having relative dimension 
two.

\begin{prop}\label{Proposition 4.8}
Let $Y$ be a $\Q$-homology threefold with a $G_m$-action. Suppose that the quotient morphism has equi-dimension two. 
Let $q : Y \to C$ be the quotient morphism, where $C$ is a smooth affine curve. Then each fiber is isomorphic to 
$\A^2$, and $C$ is isomorphic to $\A^1$. Hence $Y$ is isomorphic to $\A^3$.
\end{prop}
\Proof
Let $F$ be a fiber of $q$. Since $\dim F=2$, there is a unique fixed point $Q$ such that the closure of every orbit 
passes through $Q$. The locus $\Gamma$ of points $Q$ is the fixed point locus $Y^{G_m}$ and $q|_\Gamma : \Gamma \to 
C$ is a bijection. Hence it is an isomorphism. Since $Y$ is smooth and the local intersection multiplicity 
$i(F,\Gamma;Q)=1$, it follows that $F$ is smooth near $Q$. Hence $F$ itself is smooth. Then it is easy to show that 
$F$ is isomorphic to $\A^2$ with an elliptic $G_m$-action. By a theorem of Sathaye \cite{Sathaye}, $q : Y \to C$ is 
an $\A^2$-bundle over $C$. Since $Y$ is then contractable to $C$, the curve $C$ is $\Q$-acyclic. Hence $C \cong \A^1$ 
and $q$ is necessarily trivial. 
\QED

\begin{question}\label{Question 4.9}{\em 
Let $Y$ be a $\Q$-homology $n$-fold with a $G_m$-action. Suppose that the quotient morphism has equi-dimension $n-1$.
Is $Y$ isomorphic to the affine space $\A^n$? }
\end{question}
In fact, one can show that each fiber of the quotient morphism $q : Y \to C$ is isomorphic to $\A^{n-1}$ with the 
induced elliptic $G_m$-action. If $q$ is locally trivial in the Zariski topology, $C$ is $\Q$-acyclic. Hence 
$C \cong \A^1$ and $Y \cong \A^n$. If $Y$ is a homology $n$-fold, Question \ref{Question 4.9} has a positive answer 
which is a theorem of Kraft-Shwarz \cite[Theorem 5]{KSch}.

\begin{remark}\label{Remark 4.10}{\em 
In \cite{KR1}, a general result has been proved. Namely, let $Y$ be a homology $n$-fold. Suppose that $Y$ is dominated 
by an affine space and $Y$ is endowed with an effective action of $T=G_m^{n-2}$ such that $\dim Y^T > 0$. Then $Y$ is 
$T$-equivariantly isomorphic to the affine space $\A^n$ with a linear action of $T$. Hence a contractible threefold 
$Y$ having a non-hyperbolic $G_m$-action is $G_m$-equivariantly isomorphic to $A^3$ provided $Y$ is dominated by an 
affine space. }
\end{remark}

The following result deals with more general cases of the quotient morphism having equi-dimension one.

\begin{lem}\label{Lemma 4.11}
Let $Y$ be a smooth affine threefold with a $G_m$-action. Suppose that the quotient morphism $q : Y \to X$ has 
equi-dimension one. Then the following assertions hold.
\begin{enumerate}
\item[(1)]
If $\dim Y^{G_m}=2$, then $Y^{G_m}$ is isomorphic to $X$ and hence $Y^{G_m}$ is connected. The morphism $q$ defines a 
line bundle over $X$.
\item[(2)]
If $\dim Y^{G_m}=1$ and $e(Y)>0$, then $Y^{G_m}$ is smooth and consists of connected components $\Gamma_1, \ldots, 
\Gamma_r$, one of which is isomorphic to $\A^1$. Further, if $Y^{G_m}$ is connected, then $e(Y)=1$ and the quotient 
surface $X$ is a normal affine surface with an embedded line and at worst cyclic quotient singularities.  
\end{enumerate}
\end{lem}
\Proof
(1)\ Since $q$ does not contain a fiber component of dimension two, $Y^{G_m}$ lies horizontally to the fibration $q$. 
Hence each fiber contains a unique fixed point. This implies that each fiber is isomorphic to $\A^1$. Considering 
the tangent space $T_{Y,Q}$ and the induced tangential representation for each $Q \in Y^{G_m}$, we know that $Y^{G_m}$ 
is smooth and isomorphic to $X$. Namely, $q$ is an $\A^1$-fibration with all reduced fibers isomorphic to $\A^1$ and 
has two cross-sections $Y^{G_m}$ and a section at infinity. Hence $q$ is in fact a line bundle. 

(2)\ The morphism $q$ is then an $\A^1_*$-fibration and $Y^{G_m}$ is a smooth curve. Let $Y^{G_m}=\Gamma_1\sqcup \cdots
\sqcup \Gamma_r$ be the decomposition into connected components. Let $\ol{\Gamma}_i$ be the smooth completion of $\Gamma_i$.
Let $g_i$ be the genus of $\ol{\Gamma_i}$ and let $n_i$ be the number of points in $\ol{\Gamma}_i\setminus \Gamma_i$.
Note that $G_m$ acts on $q^{-1}(X-q(Y^{G_m}))$ without fixed points. Hence, by Lemma \ref{Lemma 4.6}, the Euler number of 
$q^{-1}(X-q(Y^{G_m}))$ is zero. Note that $q^{-1}(q(\Gamma_i))$, taken with reduced structure, is a union of two 
$\A^1$-bundles over $q(\Gamma_i)$ meeting transversally along the section $\Gamma_i$. In fact, since the fiber $q^{-1}(P)$ 
over a point $P \in q(\Gamma_i)$ is a cross with each branch meeting $\Gamma_i$ transversally in one point, 
$q^{-1}(q(\Gamma_i))$ consists of two irreducible components $W_i^{(1)}$ and $W_i^{(2)}$, each of which has an 
$\A^1$-fibration with a cross-section $\Gamma_i$. Hence $W_i^{(1)}$ and $W_2^{(2)}$ are $\A^1$ bundles over $q(\Gamma_i)$ 
meeting transversally along $\Gamma_i$. Hence the Euler number of $q^{-1}(q(\Gamma_i))$ is equal to $2-2g_i-n_i$. 
This observation yields a relation
\[
0 < e(Y)=\sum_{i=1}^r(2-2g_i-n_i).
\]
Since $\Gamma_i$ is an affine curve, we have $n_i \ge 1$. If none of $\Gamma_1, \ldots,\Gamma_r$ is isomorphic to $\A^1$,
then the right side of the above equality is less than or equal to zero, which is a contradiction. Hence one of them is 
isomorphic to $\A^1$. If $Y^{G_m}$ is connected, then $Y^{G_m}=\Gamma_1\cong \A^1$. It follows from the footnotes in the 
proof of Theorem \ref{Theorem 2.6} that $q$ induces a closed embedding of $Y^{G_m}$ into $X$ and that $X$ is smooth near 
$q(Y^{G_m})$. By Lemma \ref{Lemma 2.13}, $X$ has at worst cyclic quotient singularities in the open set 
$X\setminus(q(Y^{G_m}))$.
\QED  

\begin{remark}\label{Remark 4.12}{\em
Given a smooth affine surface $X_0$, we can produce a smooth affine surface $X$ by a half-point attachment \cite[p. 233]{Mi6} which 
contains $X_0$ as an open set and has an embedded affine line $C$. By the same procedure as in Lemma \ref{Lemma 4.1}, 
we take a product $X\times \A^1$ and apply the affine modification of $X\times \A^1$ with center $C\times (0)$ in 
$X\times (0)$. The resulting threefold is a smooth affine threefold $Y$ with a $G_m$-action such that the quotient morphism 
$q : Y \to Y\quot G_m$, where $X \cong Y\quot G_m$ and $Y^{G_m} \cong C$. Further, $\lkd(Y)=\lkd(X)$. }
\end{remark}

\mvskip

\noindent
R.V. Gurjar, School of Mathematics, Tata Institute for Fundamental Research, 400005 Homi Bhabha Road,
Mumbai, India; \\
e-mail: gurjar@math.tifr.res.in
\svskip

\noindent
M. Koras, Institute of Mathematics, Warsaw University, ul. Banacha 2, Warsaw, Poland; \\
e-mail: koras@mimuw.edu.pl
\svskip

\noindent
K. Masuda, School of Science and Technology, Kwansei Gakuin University, Hyogo 669-1337, Japan;\\
e-mail: kayo@kwansei.ac.jp
\svskip

\noindent
M. Miyanishi, Research Center for Mathematical Sciences, Kwansei Gakuin University, Hyogo 669-1337, Japan;\\
e-mail: miyanisi@kwansei.ac.jp
\svskip

\noindent
P. Russell, Department of Mathematics and Statistics, McGill University, Montreal, 805 Sherbooke St. West, Canada;\\
e-mail: russell@math.mcgill.ca

\end{document}